\documentclass[reqno]{amsart}
\usepackage{caption}
\usepackage{a4wide,amsfonts,amsmath,mathrsfs}
\usepackage{graphicx, caption, subcaption}
\usepackage{amssymb}
\usepackage{hyperref}
\usepackage{wrapfig}
\usepackage{epsfig}
\usepackage{float}
\usepackage{bbm}
\usepackage{multirow}
\usepackage[dvipsnames]{xcolor}
\usepackage{tikz}
\usetikzlibrary{decorations.pathreplacing}
\usetikzlibrary{fadings}
\usepackage[utf8]{inputenc}
\usepackage[english]{babel}
 
\long\def\MSC#1\EndMSC{\def\arg{#1}\ifx\arg\empty\relax\else
	{\narrower\noindent%
		{2010 Mathematics Subject Classification}: #1\\} \fi}
\long\def\KEY#1\EndKEY{\def\arg{#1}\ifx\arg\empty\relax\else
	{\narrower\noindent% 
		Keywords: #1\\}\fi}
		
\newcommand{\R}{\mathbb{R}}
\newcommand{\N}{\mathbb{N}}
\newcommand{\T}{{\rm T}}
\newtheorem{theorem}{Theorem}[section]
\newtheorem{remark}[theorem]{Remark}
\usepackage{multirow}
\usepackage{tikz}
\usetikzlibrary{decorations.pathreplacing}
\usetikzlibrary{fadings}

\begin{document}
\title[Neural networks for classification of stroke in 3D EIT]{Neural networks for classification of strokes in electrical impedance tomography on a 3D head model}

\author[V.~Candiani]{Valentina Candiani}
\address[V.~Candiani]{Department of Mathematics and Systems Analysis, Aalto University, P.O. Box~11100, FI-00076 Aalto, Finland.}
\email{valentina.candiani@aalto.fi}

\author[M.~Santacesaria]{Matteo Santacesaria}
\address[M.~Santacesaria]{MaLGa Center, Department of Mathematics, University of Genoa, Via Dodecaneso 35, 16146 Genova, Italy.}
\email{matteo.santacesaria@unige.it}

\begin{abstract}
We consider the problem of the detection of brain hemorrhages from three dimensional (3D) electrical impedance tomography (EIT)  measurements. This is a condition requiring urgent treatment for which EIT might provide a portable and quick diagnosis. We employ two neural network architectures -- a fully connected and a convolutional one -- for the classification of hemorrhagic and ischemic strokes. The networks are trained on a dataset with $40\,000$ samples of synthetic electrode measurements generated with the complete electrode model on realistic heads with a 3-layer structure. We consider changes in head anatomy and layers, electrode position, measurement noise and conductivity values. We then test the networks on several datasets of unseen EIT data, with more complex stroke modeling (different shapes and volumes), higher levels of noise and different amounts of electrode misplacement. On most test datasets we achieve $\geq 90\%$ average accuracy with fully connected neural networks, while the convolutional ones display an average accuracy $\geq 80\%$. Despite the use of simple neural network architectures, the results obtained are very promising and motivate the applications of EIT-based classification methods on real phantoms and ultimately on human patients.
\end{abstract}

\maketitle
\KEY
electrical impedance tomography, 
classification of brain strokes,
fully connected neural networks,
convolutional neural networks,
computational head model
\EndKEY
	
\MSC
65N21,
68T07,
35R30.
\EndMSC
	
\section{Introduction}
\label{sec:intro}
\emph{Electrical impedance tomography} (EIT) is a noninvasive imaging modality for recovering information about the electrical conductivity inside a physical body from boundary measurements \linebreak of current and potential. In practice, a set of contact electrodes is employed to drive current patterns into the object and the resulting electric potential is measured at (some of) the electrodes. The reconstruction process of EIT requires the solution of a highly nonlinear inverse problem on noisy data. This problem is typically ill-conditioned \cite{Borcea02,Cheney99,Uhlmann09} and solution algorithms need either simplifying assumptions or regularization strategies based on a priori knowledge.

In recent years, {\it machine learning} has arisen as a data-driven alternative that has shown tremendous improvements upon the ill-posedness of several inverse problems \cite{arridge2019,McCann17,Lucas18}. It has been already successfully applied in EIT imaging \cite{hamilton2018deep, hamilton2019beltrami,Li19,seo2019learning}.
The purpose of this work is to apply machine learning to the problem of classification of brain strokes from EIT data.

Stroke, a serious and acute cerebrovascular disease, and a leading cause of death, can be of two types: \textit{hemorrhagic}, caused by blood bleeding into the brain tissue through a ruptured intracranial vessel, and \textit{ischemic}, caused by vascular occlusion in the brain due to a blood clot (thrombosis). Visible symptoms are precisely the same in both cases, which makes it extremely difficult to differentiate them without advanced imaging modalities. Ischemic stroke can be treated through the use of thrombolytic (or clot-dissolving) agents within the first few hours \cite{Hacke08}, since human nervous tissue is rapidly lost as stroke progresses. On the other hand, thrombolytics are harmful, or even potentially fatal, to patients suffering from hemorrhagic stroke, thus it is essential to distinguish between the two types. A rapid and accurate diagnosis is crucial to triage patients \cite{Dowrick15, Yang17} to speed up clinical decisions concerning treatments and to improve the recovery of patients suffering from acute stroke \cite{Saver06}.

In this work we study how to accelerate the diagnosis of acute stroke using EIT. Currently, stroke can be classified only by using expensive hospital equipment such as X-ray CT. On the contrary, an EIT device is cheap, compact and could be carried in an ambulance (even though measurements would need to be taken while the patient is not moving). The main challenge, for emergency use, is that data are collected at a single time frame: this excludes time-difference imaging \cite{Barber84} and leaves absolute and frequency-difference \cite{McEwan06} imaging as the only options. Another important application, where measurements at different times are available, is bedside real-time monitoring of patients after the acute stage of stroke. In both scenarios, getting a full image reconstruction with the existing inversion algorithms is computationally heavy and time-consuming, and thus machine learning techniques can be used to expedite the process. In this work we focus on the case of absolute imaging, i.e., where EIT measurements are available at a single time frame and at a single frequency, which can be considered the most challenging scenario.

Although EIT for brain imaging has been studied for decades \cite{holder1992electrical,McEwan06,fabrizi2009electrode,malone2014stroke,nissinen2015contrast,yang2016novel,mcdermott2019bi,kolehmainen2019incorporating}, there are only few recent results that employ machine learning algorithms for stroke classification. The work \cite{McDermott18} proposes the use of both Support Vector Machines (SVM) and Neural Networks (NN) for detecting brain hemorrhages using EIT measurements in a $2$-layer model of the head. The main weakness of the model, however, is that it does not take into account the highly resistive skull layer, which is known to have a shielding effect when it comes to EIT measurements. Also, only a finite set of head shapes is considered and the model lacks the ability to generalize to new sample heads. The more recent work \cite{mcdermott2020multi} considers a 4-layer model for the heads and uses data from 18 human patients \cite{goren2018multi} that are classified using SVM. The main difference with our method is that our classification is made directly from raw electrode data, while in \cite{mcdermott2020multi} a preprocessing step involving a precise knowledge of the anatomy of the patient's head is required. Moreover, only strokes of size $20$ ml or $50$ ml are considered in four specific locations, while our datasets include strokes with volume as small as $1.5$ ml located anywhere within the brain tissue. Another methodology is shown in \cite{Agnelli20}, where first Virtual Hybrid Edge Detection (VHED) \cite{Greenleaf18} is used to extract specific features of the conductivity, then neural networks are trained to identify the stroke type. This approach is very promising but currently limited to a 2D model. Applications of deep neural networks to EIT have been also considered in \cite{Adler94,Capps2021,Lampinen99,hamilton2018deep, hamilton2019beltrami,Li19}. One could also use some machine learning techniques to form a model for the head shapes: see \cite{seo2019learning} for an approach that could potentially be applicable to head imaging.

In this work we consider two different types of NN, a fully connected and a convolutional NN, that we feed with absolute EIT measurements and produce a binary output. These measurements, which form the training and test datasets, are simulated by using the so-called \emph{complete electrode model} (CEM) \cite{Cheng89, Somersalo92} on a computational 3-layer head model, where each layer corresponds to a different head anatomical region: scalp, skull and brain tissue. 

The training and test datasets are made of pairs of simulated electrode measurements at a single time frame and a label which indicates whether the data are associated with a hemorrhagic stroke or not. More precisely, label $1$ is meant to indicate a hemorrhagic stroke, while label $0$ stands for either an ischemic stroke or no stroke.
This is motivated by the fact that detecting the presence or absence of hemorrhage may be sufficient to initiate appropriate treatment. More precisely, the datasets contain EIT measurements in the following proportions: $50\%$ hemorrhages, $25\%$ ischemic stroke and $25\%$ healthy brains. We chose to include a large number of healthy patients in order to cover a broader range of potential applications and not restrict ourselves only to the emergency setting.
The measurements in the training and test datasets are generated by varying the conductivity distribution, the electrode positions and noise, the shape of the scalp, the skull and the brain tissue. We model a hemorrhage as a volume of the brain with higher conductivity values with respect to the brain tissue, and an ischemic stroke with lower values, based on the available medical literature \cite{Latikka01, mccann2019variation}. In the training dataset the strokes are modeled as a single ball inclusion of higher or lower conductivity, while in the test datasets we consider different shapes of multiple inclusions. Concerning the variations in the geometry, a joint principal component model for the variations in the anatomy of the human head \cite{Candiani19} is considered, so that we are able to generate realistic EIT datasets for brain stroke classification on a 3D finite element (FE) head model.
The training dataset is made of $40\,000$ pairs of electrode data and labels, while every test dataset is made of approximately $5\,000$ samples. No validation was used in the training of the fully connected network, while for the convolutional one the training set was randomly split (83\% training, 17\% validation). These test datasets take into account a variety of possible errors in the measurement setup: slight variations in the background conductivity and in the contact impedances are considered, along with misplacement of electrodes and mismodeling of the head shape. The functionality of the chosen methods is demonstrated via the measures of accuracy, sensitivity and specificity of the networks trained with noisy EIT data. 

Our numerical tests show that the probability of detecting hemorrhagic strokes is reasonably high, even when the electrodes are misplaced with respect to their intended location and the geometric model for the head is inaccurate. We find that both fully connected neural networks and convolutional neural networks are efficient tools for the described classification. More precisely, in our experiments we observe that a shallow fully connected neural network generalizes better to the test datasets than a convolutional one.

This paper is organized as follows. In Section \ref{sec:datagen} we recall the CEM and the parametrized head model with the workflow for mesh generation. Neural networks and their specifications are introduced in Section \ref{sec:NN}. Section \ref{sec:settings} presents the experiment settings, while numerical results are described in Section \ref{sec:results}. Finally, Section \ref{sec:conclusion} lists the concluding remarks.

\section{Generation of electrode data}
\label{sec:datagen}
	
\subsection{Forward model}
We start by recalling the CEM~\cite{Somersalo92} of EIT. Let $\Omega \subset \R^3$ denote a bounded Lipschitz domain and assume there is a set of $M \in \N \setminus \{ 1 \}$ contact electrodes $E_1, \dots, E_M$ attached to its boundary $\partial\Omega$. When a single measurement by an EIT device is performed, net currents $I_m\in\R $, $m=1, \dots, M$ are driven through each $E_m$ and the resulting constant electrode potentials $U_m \in \R$, $m=1, \dots, M$, are measured. Due to conservation of electric charge, any applicable current pattern $I = [I_1,\dots,I_M]^{\rm T}$ belongs to the mean-free subspace
\[ 
\R^M_\diamond \, := \, \Big\{J \in\R^M\,\Big|\, \sum_{m=1}^M J_m = 0\Big\}.
\]

The contact impedances at the electrode-object interfaces are modeled by a vector \linebreak $z = [z_1, \dots, z_M]^{\rm T}\in \R_+^M$. The electrode patches are identified with the nonempty, connected and open subsets of $\partial\Omega$ that they cover and assumed to be well separated, i.e.,~$\overline{E}_m\cap \overline{E}_l = \emptyset$ if $m\not= l$. We denote $E = \cup E_m$. 
The electromagnetic potential $u$ inside $\Omega$ and the piecewise constant potentials on the electrodes $U$ weakly satisfy

\begin{equation}
\label{eq:cemeqs}
\begin{array}{ll}
\displaystyle{\nabla \cdot(\sigma\nabla u) = 0 \qquad}  &{\rm in}\;\; \Omega, \\[6pt] 
{\displaystyle {\nu\cdot\sigma\nabla u} = 0 } \qquad &{\rm on}\;\; \partial \Omega\setminus \overline{E}, \\[6pt]{\displaystyle {u+z_m\nu\cdot\sigma\nabla u} = U_m } \qquad &{\rm on}\;\; E_m,  \qquad m=1,\ldots,M, \\[2pt] 
{\displaystyle \int_{E_m}\nu\cdot\sigma\nabla u\,{\rm d}S} = I_m, \qquad & m=1,\ldots,M, \\[4pt]
\end{array}
\end{equation}
where $\nu \in L^\infty(\partial \Omega, \R^3)$ denotes the exterior unit normal of $\partial \Omega$.  Moreover, the isotropic conductivity distribution $\sigma$ describing the electric properties of $\Omega$ is assumed to belong to 
\begin{equation}
\label{eq:sigma}
L^\infty_+(\Omega) := \{ \varsigma \in L^\infty(\Omega) \ | \ {\rm ess} \inf \varsigma > 0 \}.
\end{equation}

A physical justification of \eqref{eq:cemeqs} can be found in \cite{Cheng89}. Given an input current pattern $I \in \R_\diamond^M$, a conductivity $\sigma$ and contact impedances $z$ with the properties described above, the pair $(u, U) \in H^1(\Omega)  \oplus \R_{\diamond}^M$ is the unique solution of the elliptic boundary value problem \eqref{eq:cemeqs} according to \cite{Cheng89, Somersalo92}. Note that the use of $\R_{\diamond}^M$ corresponds to systematically choosing the ground level of potential so that the mean of the electrode potentials is zero. 
The measurement, or current-to-voltage map of CEM is defined as the mapping $ I \mapsto U$, from $\R_\diamond^M$ to $\R_\diamond^M$.

\subsection{Head model}
\label{ssec:head}
\begin{figure}[H]
 \center{
  {\includegraphics[width=5.5cm]{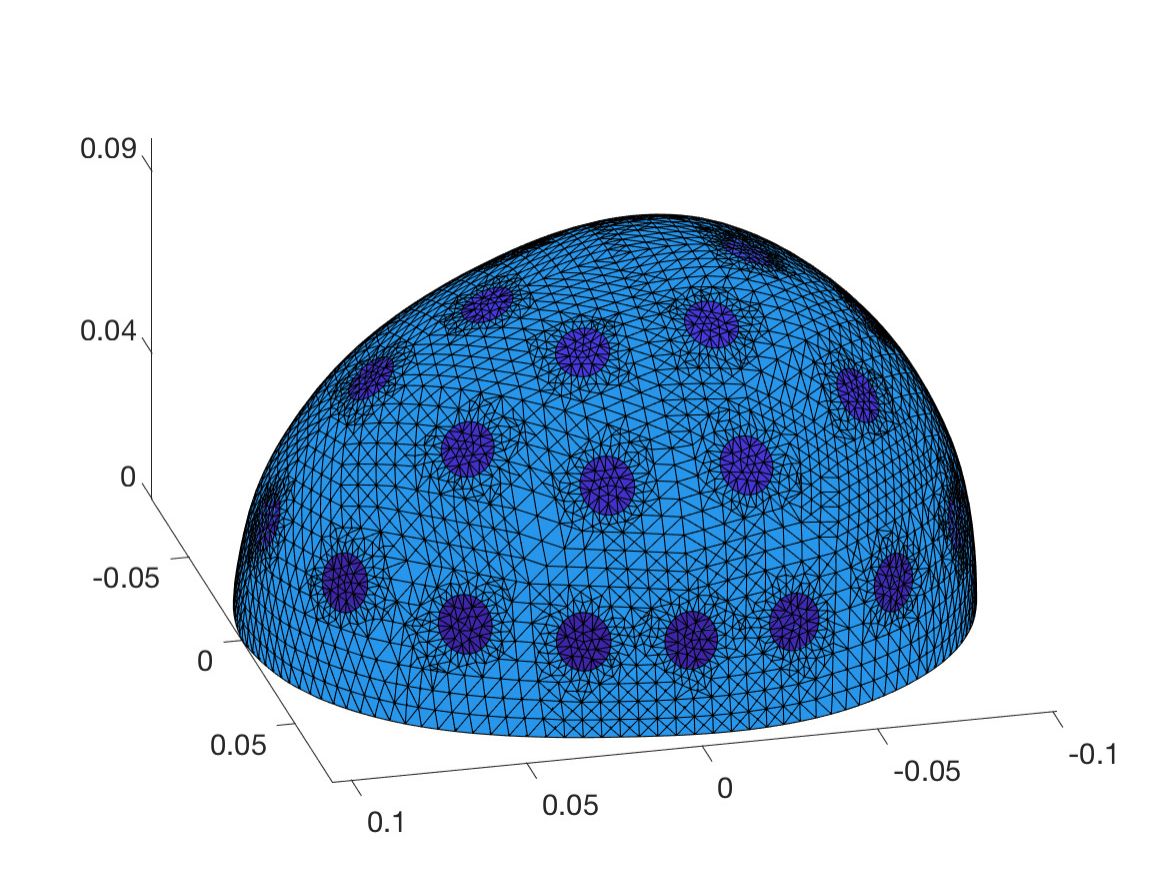}}
  \quad
  {\includegraphics[width=5.5cm]{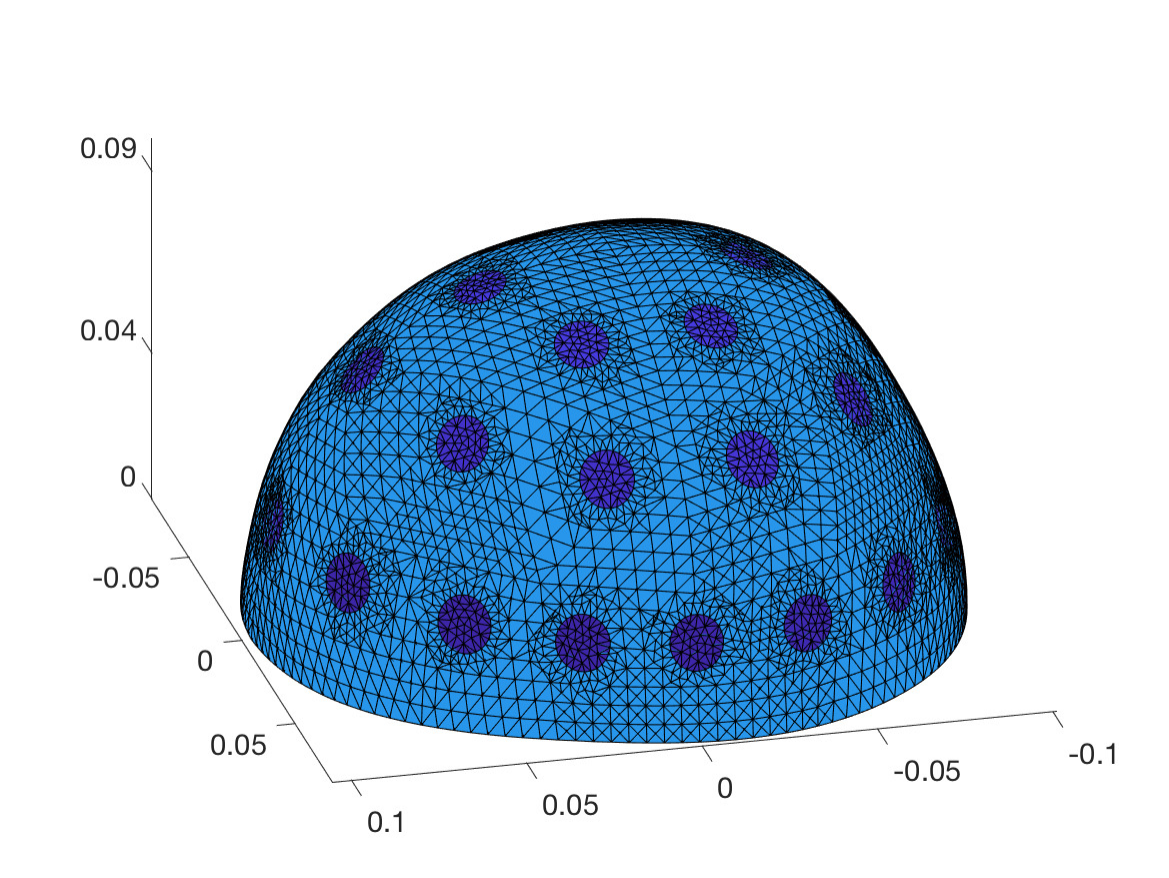}}
  }
  \center{
  {\includegraphics[width=5.2cm]{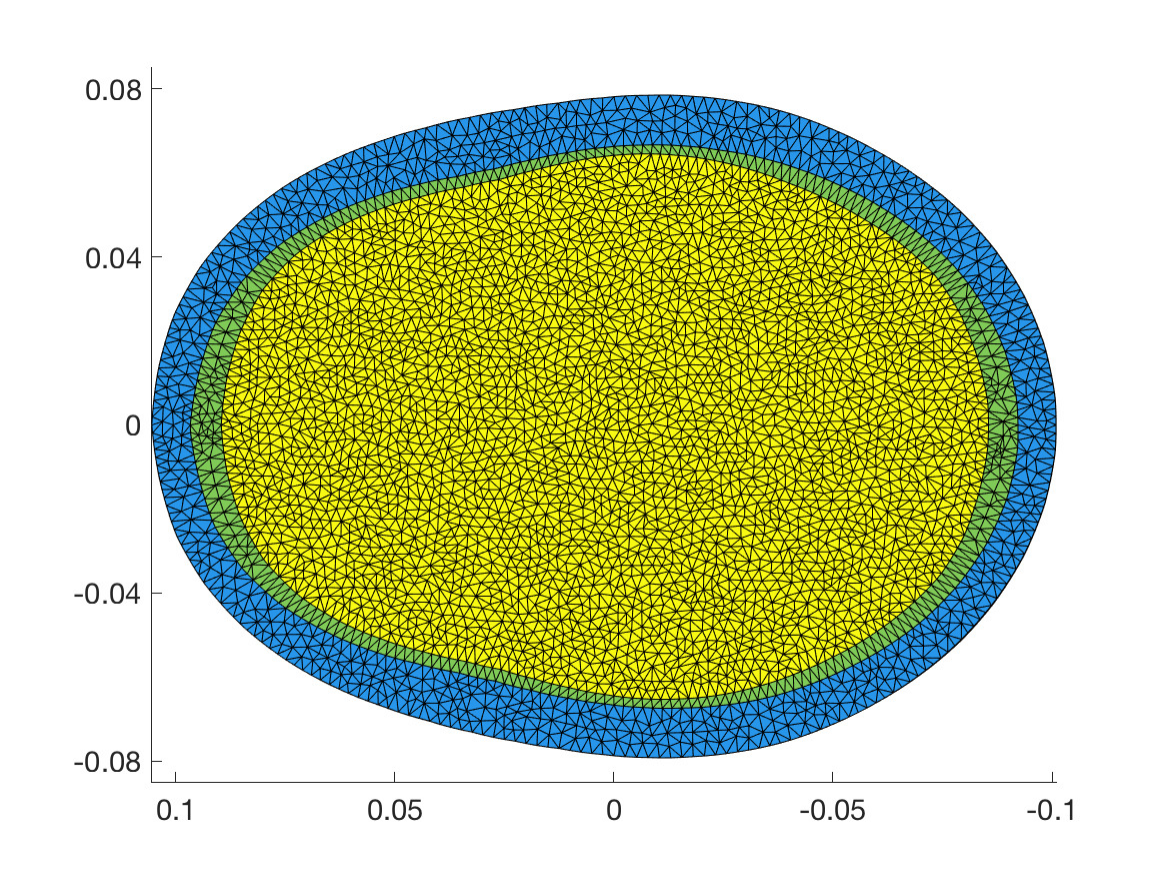}}
  \quad
  {\includegraphics[width=5.2cm]{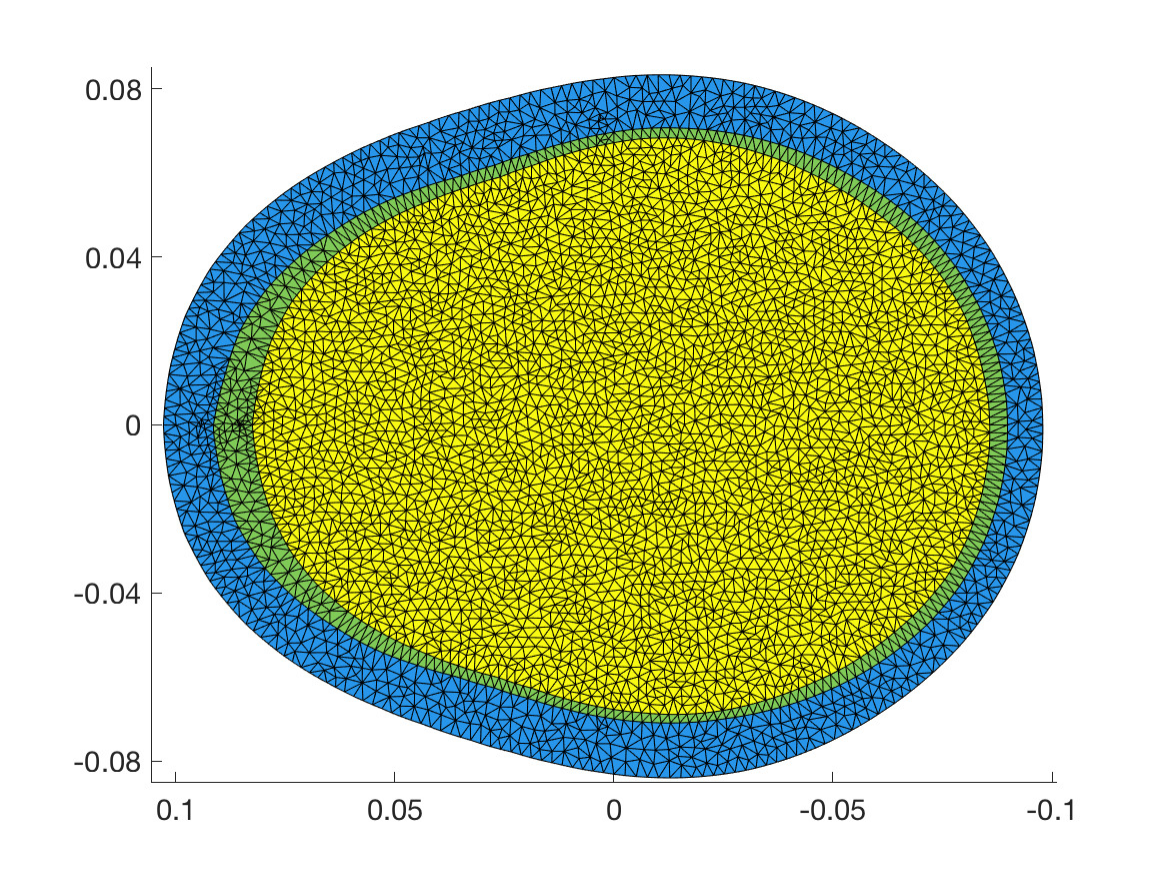}}
  }
 \caption{Top row: two different head models, oriented with the forehead on the left and the back of the head on the right. The $32$ electrodes are at their intended positions, with the FE mesh associated to the head appropriately refined around them. Bottom row: the corresponding bottom face of each head model, where the three layers associated to scalp, skull and brain tissues are visible. Note that the first model has a more flattened forehead (top row) and it is narrower in the coronal plane (bottom row, $y$ direction), corresponding to a shorter distance between the ears. Also, the thickness of the scalp layer is clearly different. The unit of length is meter.}
\label{fig:head}
\end{figure}

The head model used in this work follows the same approach as in \cite{Candiani19}, though slightly modifying and upgrading the setting to a three-layers model. We define a {\em layer} for each one of the anatomical structures that we are considering for this particular head model. There are $L = 3$ different layers: the scalp layer, i.e., the outer one corresponding to the skin, the resistive skull layer and the interior brain layer (see Figure \ref{fig:head}).

For each layer, the library of $n=50$ heads from \cite{Lee16} is used to build the model for the variations in the shape and size of the human head. We can represent the crown of the $l$th layer in the $j$th head, for $l=1,2,3$ and $j=1,\dots,n$, as the graph of a function
\begin{equation}
  \label{eq:jth_head}
S_j^l:
\left\{
\begin{array}{l}
\mathbb{S}_+ \to \R^3, \\[1mm]
\hat{x} \mapsto r_j^l(\hat{x}) \, \hat{x},
\end{array}
\right.
\end{equation}
where $\mathbb{S}_+$ is the upper unit hemisphere, i.e.,
$$
\mathbb{S}_+ = \big\{ x \in \R^3 \; | \; \| x \|_2 = 1 \ {\rm and} \ x_3 > 0 \big\},
$$
and $r_j^l:  \mathbb{S}_+ \to \R_+$ gives the distance from the origin to the surface of the $l$th layer of the $j$th head as a function of the direction $\hat{x} \in \mathbb{S}_+$, where the origin is set at approximately the center of mass of each bottom face of the heads (see Figure \ref{fig:head}). 

Then, for each layer $l$, we introduce the average head and perturbations
\begin{equation}
  \label{eq:mean_etc}
\bar{r}^l = \frac{1}{n} \sum_{j=1}^n r_j^l \qquad {\rm and} \qquad \rho_j^l = r_j^l - \bar{r}^l, \quad j=1, \dots, n, \,\  l = 1, 2, 3,
\end{equation}
where $\bar{r}^l$ describes the $l$th layer of the average head and $\rho_{1}^l, \dots, \rho_n^l$ are the corresponding perturbations that define the employed library of heads. We assume the functions $\rho_1^l, \dots, \rho_n^l $ belong to $H^1(\mathbb{S}_+)$ and are linearly independent for every $l$.

Mimicking the formulation in \cite{Candiani19}, we introduce a joint principal component model involving all the three layers. A single head in the library defines an object in the space $[H^1(\mathbb{S}_+)]^3$, that is, a three dimensional vector whose components are in the space $H^1(\mathbb{S}_+)$ and define the three layers. The reason for choosing $H^1(\mathbb{S}_+)$ is two-fold. First of all, according to the numerical tests in \cite{Candiani19}, the simplest option $s=0$ in $H^s(\mathbb{S}_+)$ leads to undesired cusps in the optimal basis $\hat{\rho}_1, \dots, \hat{\rho}_n$ defined below. On the other hand, $s\in \N\setminus \{1\}$ would require the use of higher order conforming finite elements and the implementation of the needed inner products on $\mathbb{S}_+$ for $s\in\R_+\setminus \N$ would lead to unnecessary technical considerations. The construction of the new principal component model for the head library is performed as in \cite{Candiani19}, with the exception that the inner product between two elements of $[H^1(\mathbb{S}_+)]^3$, say, $v$ and $w$, is defined as
\begin{equation}
\label{eq:inner}
(v, w)_{[H^1(\mathbb{S}_+)]^3} := \sum_{i=1}^{3}(v_i, w_i)_{H^1(\mathbb{S}_+)}.
\end{equation}

Our aim is to look for an $\tilde{n}$-dimensional subspace $V_{\tilde{n}} \subset  [H^1(\mathbb{S}_+)]^3$, $1 \leq \tilde{n} \leq n$,
that satisfies
\begin{equation}
  \label{eq:minim}
\sum_{j=1}^{n} \min_{\eta \in V_{\tilde{n}}} \| \rho_j - \eta \|_{[H^1(\mathbb{S}_+)]^3}^2 \leq \sum_{j=1}^{n} \min_{\eta \in W} \| \rho_j - \eta \|_{[H^1(\mathbb{S}_+)]^3}^2,
\end{equation}
for all $\tilde{n}$-dimensional subspaces $W$ of the Sobolev space $[H^1(\mathbb{S}_+)]^3$. The purpose is to find a low dimensional subspace that on average contains the best approximations for the perturbations $ \{ \rho_j = (\rho^1_j ,\rho^2_j ,\rho^2_j ) \}_{j=1}^n$, where the quality of the fit is measured by the squared norm of $[H^1(\mathbb{S}_+)]^3$.

Following the approach in \cite{Candiani19}, we define the matrix $R$ that takes into account the variations in every layer: 
$$
R_{ij} = \sum_{l=1}^{3}( \rho_i^l, \rho_j^l)_{H^1(\mathbb{S}_+)} = (\rho_i, \rho_j)_{[H^1(\mathbb{S}_+)]^3}, \qquad i,j = 1, \dots, n.
$$
By applying Lemma 3.1 from \cite{Candiani19} with minor modifications, we obtain the following set of orthonormal basis functions for $V_{\tilde{n}}$:
\begin{equation}
\label{eq:rhohat}
\hat{\rho}_k:=(\hat{\rho}_k^1, \hat{\rho}_k^2, \hat{\rho}_k^3) \quad \text{with} \quad \hat{\rho}_k^l := w_k^{\rm T}\rho^l \qquad k=1, \dots, \tilde{n},
\end{equation}
where 
\begin{equation}
\label{eq:sceig}
w_k = \frac{1}{\sqrt{\lambda_k}} \, v_k, \qquad k=1, \dots, \tilde{n},
\end{equation}
$\lambda_k, v_k$ are eigenvalues and orthonormal eigenvectors of $R \in \R^{n \times n}$ and we have defined $\rho^l = [\rho_1^l, \dots, \rho_n^l]^{\rm T}: \mathbb{S}_+ \to \R^n$, for $l = 1, 2, 3$. The positive eigenvalues $\lambda_k$ are listed in descending order, and the corresponding $l$-dependent eigenvectors are employed in the definition of $\hat{\rho}_k^l$ for all $l = 1, 2, 3$.

The parametrization for the $l$th layer in our head model can then be written as
\begin{equation}
\label{eq:theparam_l}
S^l(\hat{x}; \alpha) = \Big( \bar{r}^l(\hat{x}) +
\sum_{k=1}^{\tilde{n}} \alpha_k \hat{\rho}_k^l(\hat{x}) \Big) \hat{x}, \qquad \hat{x} \in \mathbb{S}_+, \,\ l = 1, 2, 3,
\end{equation}
where $\hat{\rho}_k^l$ are defined as in \eqref{eq:rhohat}, $\alpha_k$ are free shape coefficients and $1\leq\tilde{n}\leq n$ is chosen appropriately (cf.~\cite{Candiani19}). When generating random head structures for one numerical experiment, the vector of shape coefficients $\alpha \in \mathbb{R}^{\tilde n}$ is drawn from a Gaussian distribution $\mathcal{N}(0, \Gamma_\alpha)$, with the covariance matrix: 
\begin{equation}
  \label{eq:Gamma_alpha3}
(\Gamma_\alpha)_{kh} = \frac{1}{\tilde n -1}\sum_{j=1}^{\tilde n}(\rho_j, \hat{\rho}_k)_{[H^1(\mathbb{S}_+)]^3} (\rho_j, \hat{\rho}_h)_{[H^1(\mathbb{S}_+)]^3};
\end{equation}
see once again \cite{Candiani19} for further details.

\subsection{Generating the FEM mesh}
Our workflow for generating a tetrahedral mesh for the head model consists of three steps: generation of an initial surface mesh, insertion of electrodes and tetrahedral mesh generation. The initial surface mesh is constructed by subdividing $k$ times a coarse surface partition consisting  of four triangles, where $k \in \N$ can be chosen by the operator of the algorithm, then $M$ electrodes are inserted in the surface mesh following the process described in \cite{Candiani19}. A (dense) mesh $\mathcal{T}_m$ for the resulting polygonal domain is generated using the Triangle software \cite{Shewchuk96}. After inserting all the electrodes, the process is completed by generating a tetrahedral partition for the whole volume by TetGen \cite{Hang15} starting from the formed surface mesh. 

With the current head library, the average head size is approximately $20\times16$ cm in the axial plane, while its height is about $9$ cm. The head shapes obtained by changing the shape parameters in \eqref{eq:theparam_l} are variations from this average, up to a maximum of approximately $2$ cm difference in each dimension. The thickness of the scalp varies within the range $10-20$ mm, while the skull is about $2-15$ mm thick. Throughout all the experiments, the number of electrodes is chosen to be $M=32$ and we select $\tilde n = 10$.

\begin{remark}
Despite being an upgrade with respect to the computational head model used in \cite{Candiani19}, this three-layer  model is clearly still a simplified version of the true head anatomy. In particular, it does not take into account the shunting effect of the highly conductive cerebrospinal fluid (CSF) layer inside the skull. The CSF is known to represent a major challenge in EIT brain imaging, for it is extremely difficult to distinguish it from a bleed.
\end{remark}

\section{Neural networks}
\label{sec:NN}

There is a wide range of machine learning classification algorithms available in the literature. We chose to use neural networks (NN) since they consistently outperformed kernel methods for this specific kind of nonlinear dataset in our preliminary numerical tests.

We consider a 2-layer fully connected network and a convolutional neural network for the classification of brain strokes from electrode data. In both cases, as detailed in Section~\ref{sec:settings}, the input is a vector of size $M(M-1)$, which represents a single set of electrode measurements, where $M$ is the number of electrodes. The output layer is a single scalar value between $0$ and $1$. A rounding is then applied to the output in order to obtain a binary value for the classification. Since in our experiments $M = 32$, the input layer is composed of $992$ neurons.

\subsection{Fully connected neural networks}\label{sub:FCNN}
We consider a fully connected neural network (FCNN) that takes as input electrode data generated as discussed in Section~\ref{sec:datagen} and gives a binary output: $1$ for hemorrhage, $0$ for no hemorrhage.

\begin{figure}[H] 
\begin{center}
	\begin{tikzpicture}[shorten >=1pt]
		\tikzstyle{unit}=[draw,shape=circle,minimum size=.5cm]
		%\tikzstyle{hidden}=[draw,shape=circle,fill=black!25,minimum size=1.15cm]
		\tikzstyle{hidden}=[draw,shape=circle,minimum size=.5cm]
 
		\node[unit](x0) at (3,3){$ $};
		\node[unit](x1) at (3,2.2){$ $};
		\node at (3,1.6){\vdots};
		\node[unit](xd) at (3,0.8){$ $};
 
%		\node[hidden](h10) at (5,4){$ $};
		\node[hidden](h11) at (5,2.3){$ $};
		\node at (5,1.8){\vdots};
		\node[hidden](h1m) at (5,1.1){$ $};

		\node[unit](y1) at (7,1.7){$ $};

		\draw[->] (x0) -- (h11);
		\draw[->] (x0) -- (h1m);
 
		\draw[->] (x1) -- (h11);
		\draw[->] (x1) -- (h1m);
 
		\draw[->] (xd) -- (h11);
		\draw[->] (xd) -- (h1m);
 
		\draw[->] (h11) -- (y1);
		\draw[->] (h1m) -- (y1);

		\draw [decorate,decoration={brace,amplitude=10pt},xshift=-4pt,yshift=0pt] (2.8,3.3) -- (3.45,3.3) node [black,midway,yshift=+0.6cm]{input layer (992)};
		\draw [decorate,decoration={brace,amplitude=10pt},xshift=-4pt,yshift=0pt] (4.8,2.6) -- (5.45,2.6) node [black,midway,yshift=+0.6cm]{hidden layer (7)};
		\draw [decorate,decoration={brace,amplitude=10pt},xshift=-4pt,yshift=0pt] (6.8,2) -- (7.45,2) node [black,midway,yshift=+0.6cm]{output (1)};
	\end{tikzpicture}

\caption{An illustration of the architecture of our FCNN, with one input layer with $992$ nodes, one hidden layer with $7$ nodes and one output layer with a single node for the binary classification.}
\label{fig:NN}
\end{center}
\end{figure}
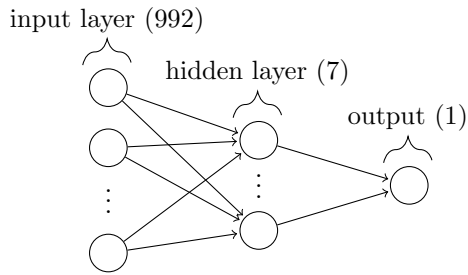

Our FCNN has two layers with weights, as shown in Figure~\ref{fig:NN}. The input layer has $992$ neurons, while the second and the third layer have $7$ and $1$ neurons, respectively. The size of the hidden layer has been chosen on the basis of the results obtained in preliminary tests on smaller datasets.
The network can be represented by the following real-valued function:
\begin{equation}\label{def:fcnn}
f(x,\theta) = g(W_2(g(W_1 x+b_1))+b_2),
\end{equation}
where we denoted by $\theta = \{W_1,W_2,b_1,b_2\}$ the set of weights and biases and
\begin{itemize}
\item $x \in \R^{992}$ is the input electrode data,
\item $W_1 : \R^{992} \to \R^{7}$, $b_1 \in \R^{7}$ are the weights and the bias of the first layer,
\item $W_2 : \R^{7} \to \R$, $b_2 \in \R$ are the weights and the bias of the second layer,
\item $g: \R \to \R$ is the sigmoid function $g(t) = 1/(1+e^{-t})$ applied component-wise.
\end{itemize}

The network $f(x,\theta)$ is trained by minimizing the cross-entropy loss
\begin{equation}\label{def:BCE}
L(\theta) = -\sum_{j}\left(y_j \log (f(x_j,\theta))+(1-y_j) \log (1-f(x_j,\theta))\right),
\end{equation}
where the sum is over the samples $x_j$ in the training set, with $y_j$ being the corresponding true binary label.

\subsection{Convolutional neural networks}

We also consider a convolutional neural network (CNN) for our classification task. We refer to \cite{lecun2010convolutional} for more details on the architecture of a CNN.

\begin{figure}[H]  
\begin{center}
\begin{picture}(300,180)
\put(0,0){\includegraphics[width=11cm]{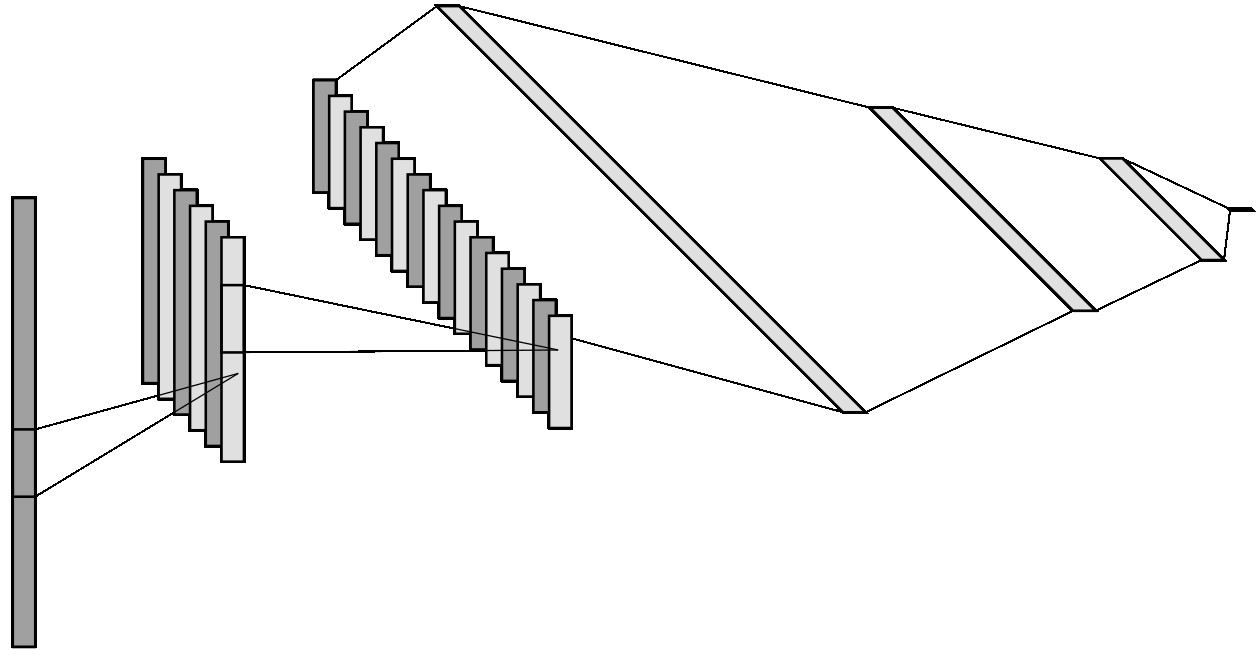}}
\put(12,100){\tiny $992$ }
\put(10,51){\tiny $3\times 1$}
\put(63,80){\tiny $3\times 1$}
\put(50,117){\tiny $6$}
\put(112,117){\tiny $16$}
\put(143,137){\tiny $240$}
\put(245,117){\tiny $120$}
\put(290,120){\tiny $84$}
\put(315,109){\tiny $1$}
\put(15,30){\tiny Max pooling}
\put(70,60){\tiny Max pooling}
\put(157,60){\tiny ReLU}
\put(237,62){\tiny ReLU}
\put(287,82){\tiny ReLU}
\end{picture}
\caption{A visual illustration of the architecture of our CNN, with two convolutional layers and four fully connected layers.}
\label{fig:CNN}
\end{center}
\end{figure}
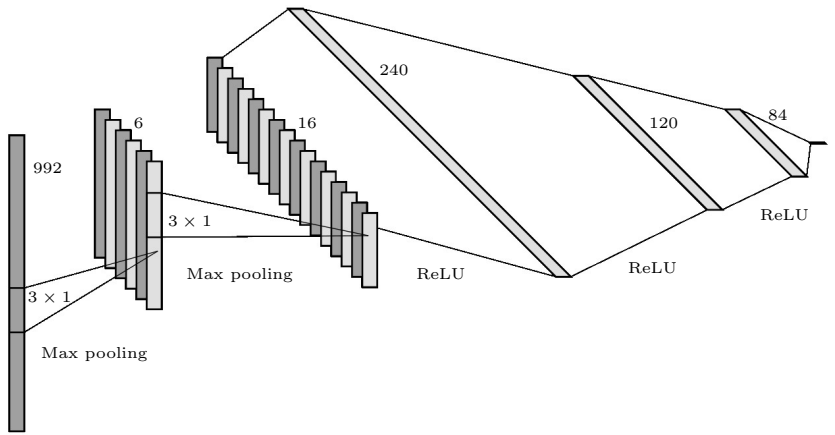

As depicted in Figure~\ref{fig:CNN}, our network has six layers with weights. The first two are convolutional and the last four are fully-connected. As with the FCNN, our CNN is trained by minimizing the cross-entropy loss \eqref{def:BCE}. The architecture was motivated by similar CNNs used in image classification \cite{lecun2010convolutional}, and the hyperparameters have been chosen after preliminary tests on smaller datasets.

Regarding the two convolutional layers, we chose to use 1D kernels. This choice might not be optimal since we are losing some geometric information about the electrode configuration. On the other hand, even though a single set of electrode measurements can be represented as a matrix of size $32\times 31$, there is no obvious advantage in considering it as an image.

In the first layer of our CNN, the $992 \times 1$ input vector is filtered with $6$ kernels of size $3 \times 1$ with stride $1$ and zero padding. Then a max-pooling layer, with kernel and stride of size $2$ and zero padding, is applied to each output channel of the first layer. This is then filtered, in the second convolutional layer, with $16$ kernels of size $3 \times 1$ (with stride $1$ and zero padding). Then another max-pooling layer, with kernel and stride of size $2$ and zero padding, is applied to the output. The third, fourth, fifth and sixth layers are fully connected and have sizes $240$, $120$, $84$ and $1$, respectively (see Section~\ref{sub:FCNN} for more details). The rectified linear unit (ReLU) activation function $g(x) = \max(0,x)$, is applied to the output of each fully connected layer, except the last one.

\section{Experimental settings}
\label{sec:settings}

\subsection{Simulation of measurement data}
\begin{figure}[H]
 \center{
  {\includegraphics[width=6.2cm]{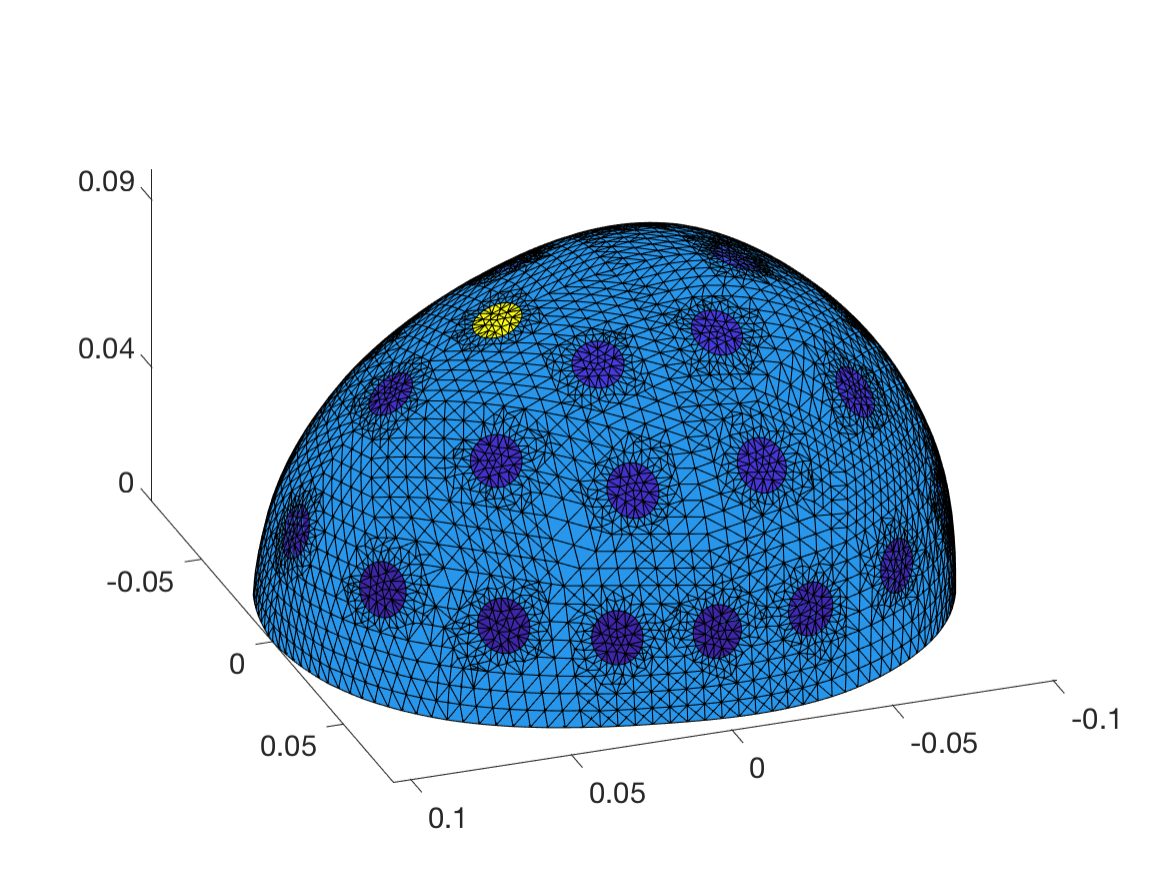}}
  \qquad
  {\includegraphics[width=6.2cm]{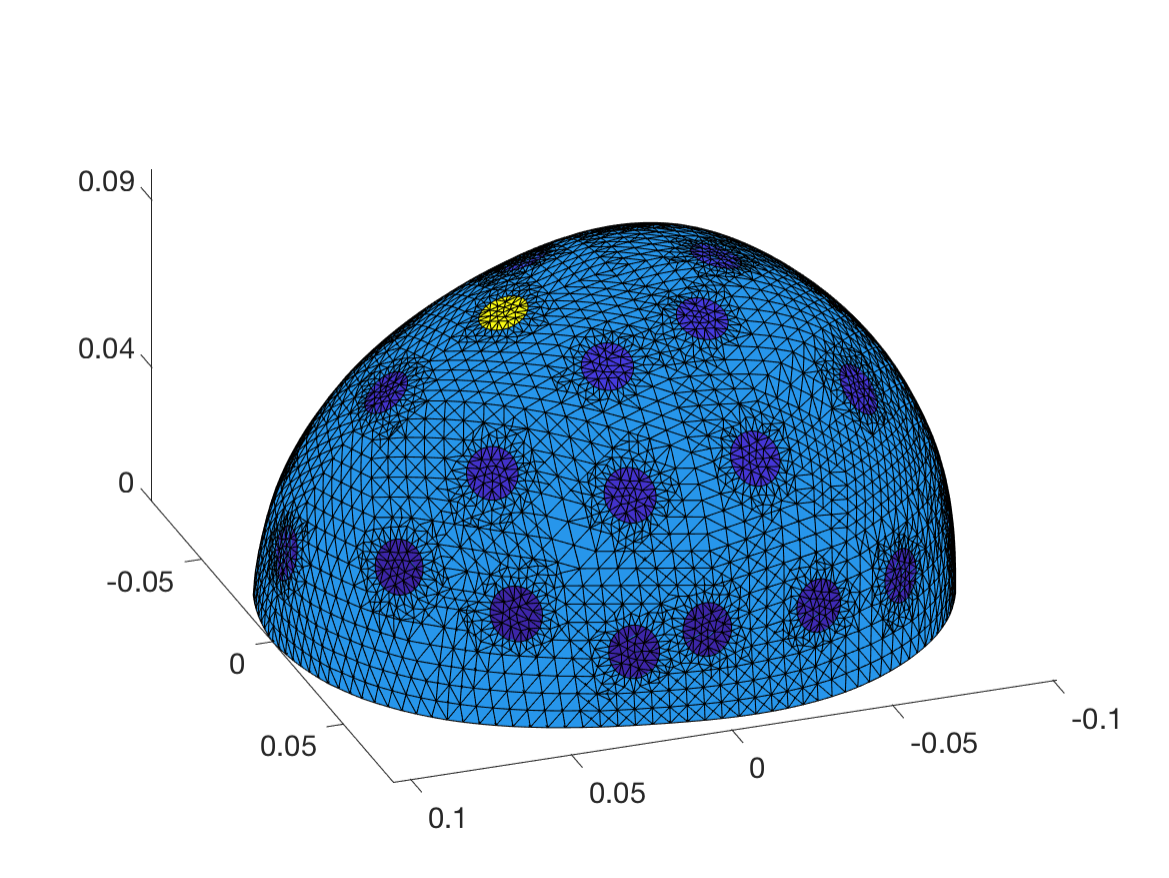}}
  }
\caption{Left: a head model with electrodes at their intended positions. Right: the same head model with misplaced electrodes, where $\varsigma_\theta = \varsigma_\phi = 0.03$ in \eqref{eq:Gamma_thph}. Counting upwards from the bottom belt, there are altogether $M =16 + 10 + 6 = 32$ electrodes of radius $R = 0.75$ cm. The current-feeding electrode is $p = 27$, i.e., the frontal one on the top belt of electrodes, highlighted in yellow. The FE mesh is refined appropriately around the electrodes. The unit of length is meter.}
\label{fig:electrodes}
\end{figure}

We assume to be able to drive $M-1$ linearly independent current patterns $I^{(1)}, \dots, I^{(M-1)} \in \R_{\diamond}^M$ through the $M$ electrodes and measure the corresponding noisy electrode potentials $V^{(1)}, \dots, V^{(M-1)} \in \R^{M}$. The employed current patterns are of the form $I^{(j)} = e_p-e_j$, $j=1, \dots, p-1, p+1, \dots, M$, where $p \in \{1, \dots, M\}$ is the label of the so-called current-feeding electrode and $e_k$ denotes a standard basis vector. Such current patterns have been used in \cite{Candiani19} and with real-world data in \cite{Darde13a} and \cite{Kourunen}. In our tests, the current-feeding electrode is always the frontal one on the top belt of electrodes (cf. Figure \ref{fig:electrodes}). The potential measurements are stacked into a single vector
\begin{equation}
\label{data}
\mathcal{V} := \big[(V^{(1)})^{\rm T},\ldots, (V^{(M-1)})^{\rm T}\big]^\T\in \R^{M(M-1)},
\end{equation}
and we  analogously  introduce the stacked forward map
$$
\mathcal{U}: \R_+^N \times \R_+^M \times \R^{\tilde{n}} \times (0,\pi/2)^M \times [0,2\pi)^M \to \R^{M(M-1)}
$$
via
$$
\mathcal{U}(\sigma, z, \alpha, \theta, \phi)
= \big[U(\sigma, z, \alpha, \theta, \phi; I^{(1)})^{\rm T},\ldots, U(\sigma, z, \alpha, \theta, \phi; I^{(M-1)})^{\rm T}\big]^\T.
$$

Here, the conductivity $\sigma \in \R_+^N$ is identified with the $N \in \N$ degrees of freedom used to parametrize it, i.e.,~the number of nodes in the mesh, the contact impedances are identified with the vector $z = [z_1, \dots, z_M] \in \R_+^M$, $\alpha \in \R^{\tilde{n}}$ is the parameter vector in \eqref{eq:theparam_l} determining the shape of the computational head model, and $\theta \in (0,\pi/2)^M$ and  $\phi \in [0,2\pi)^M$ define the polar and azimuthal  angles of the electrode center points, respectively. 

For each forward measurement, the parameters $\alpha \in \R^{\tilde{n}}$ defining the shape of the head are drawn from the distribution $\mathcal{N}(0, \Gamma_\alpha)$, where $\Gamma_\alpha \in \R^{\tilde{n} \times \tilde{n}}$ is the diagonal covariance matrix defined componentwise by \eqref{eq:Gamma_alpha3}; for a motivation of this choice, as well as for a proof of $\Gamma_\alpha$ being diagonal, see the principal component construction in \cite[Section 3.2]{Candiani19}.

Every set of measurements is performed with $M=32$ electrodes of radius $R=0.75$ cm, organized in three belts around the head (see Figure~\ref{fig:electrodes}). The expected values for the polar and azimuthal angles of the electrode centers, $\bar{\theta}$ and $\bar{\phi}$, correspond to the correct angular positions of the electrodes, i.e.,~ the positions where one originally aims to attach the electrode patches. The actual central angles of the target electrodes, i.e.,~$\theta$ and $\phi$, are then drawn from the distributions $\mathcal{N}(\bar{\theta}, \Gamma_\theta)$ and $\mathcal{N}(\bar{\phi}, \Gamma_\phi)$, where
\begin{equation}
\label{eq:Gamma_thph}
\Gamma_\theta = \varsigma_\theta^2 \mathbb{I} \qquad {\rm and} \qquad
\Gamma_\phi = \varsigma_\phi^2 \mathbb{I}. 
\end{equation}
Here, $\mathbb{I} \in \R^{M \times M}$ is the identity matrix and $\varsigma_\theta,\varsigma_\phi>0$ determine the standard deviations in the two angular directions. Notice that $\varsigma_\theta$ and $\varsigma_\phi$ must be chosen so small that the electrodes are not at a risk to overlap or move outside the crown of the computational head. 

The relative contact impedances $z_m \in \R_+$, $m=1, \dots, M$ are independently drawn from $\mathcal{N}(\bar{z}, \varsigma_z^2)$, where $\bar{z}> 0$ is chosen so much larger than $\varsigma_z>0$ that negative contact impedances never occur in practice.

Finally, the conductivity distributions are defined as follows. For each layer (scalp, skull and brain) and stroke type (hemorrhagic or ischemic) we draw parameters $\sigma_{scalp},\sigma_{skull},\sigma_{brain}$ and $\sigma_{h}$ or $\sigma_{i}$ from Gaussian distributions of the form $\mathcal{N}(\bar{\sigma}, \varsigma_\sigma^2)$ (see Table~\ref{tab:values}). With these parameters we construct a conductivity that is constant on each layer and in the stroke region (if present). The final conductivity $\sigma \in \R_+^{N}$ corresponds to a piecewise linear parametrization on a dense FE mesh with $N \approx 20\,000$ nodes and about $85\,000$ tetrahedra associated to the head defined by the shape parameters $\alpha$ and refined appropriately around the electrodes, whose positions are determined by $\theta$ and $\phi$. The number of principal components for the shape parameters $\alpha_1, \dots, \alpha_{\tilde{n}}$ (see Section~\ref{ssec:head}) is chosen to be $\tilde{n}=10$ in all the experiments.

\begin{table}[H]
	\begin{center}
		\caption{Parameter values related to the training set. These same values are also used in the simulation of the test sets if not explicitly mentioned otherwise.}
		\begin{tabular}{ll}
			\hline
			\multicolumn{2}{c}{Parameter values for the training set}\\ \hline
			conductivity of the scalp $\sigma_{\rm scalp}$ & $\sim \mathcal{N}(0.4, 0.006) $ (S/m) \\ %\hline
			conductivity of the skull $\sigma_{\rm skull} $ &  $\sim \mathcal{N}(0.06, 0.004) $  (S/m)  \\ %\hline
			conductivity of the brain $\sigma_{\rm brain}$ &  $ \sim \mathcal{N}(0.2, 0.016) $ (S/m)  \\ %\hline
			conductivity of hemorrhagic stroke $\sigma_{\rm h}$ & $ \sim \mathcal{N}(1.8, 0.012) $ (S/m) \\ %\hline
			conductivity of ischemic stroke $\sigma_{\rm i} $ &  $\sim \mathcal{N}(0.12, 0.008) $ (S/m)  \\ %\hline
			radius of the inclusion $r$ & $\sim \mathcal{U}(0.7, 2.3) $ (cm) \\ %\hline
			height of the cylindrical inclusion $h$ & $\sim \mathcal{U}(1, 3) $ (cm) \\ %\hline
			standard deviation for electrode positions $\varsigma_\theta, \varsigma_\phi $ & $\in \{ 0.01, 0.02, 0.03\} $ (rad)\\ %\hline
			contact impedances values $z_m$ & $\sim \mathcal{N}(0.01, 0.001) $ (S/m$^2$)  \\ %\hline
			shape parameters $\alpha$ & $ \sim \mathcal{N}(0, \Gamma_\alpha) $ \\ %\hline
			relative noise level $\varsigma_\eta$ & $10^{-3} $ \\ %\hline
			stroke location & uniformly random within the brain tissue \\ \hline 
		\end{tabular}
		\label{tab:values}
	\end{center}
\end{table}

We approximate the ideal data $\mathcal{U}(\sigma, z, \alpha, \theta, \phi)$ by FEM with piecewise linear basis functions and denote the resulting (almost) noiseless data by $\mathcal{U} \in \R^{M(M-1)}$. The actual noisy data is then formed as
\begin{equation}
\label{eq:noisydata} 
\mathcal{V} = \mathcal{U} + \eta,
\end{equation}
where $\eta \in \R^{M(M-1)}$ is a realization of a zero-mean Gaussian with the diagonal covariance matrix
\begin{equation}
\label{eq:noisecov} 
\Gamma_\eta = \varsigma_\eta^2 \big( \max_{j} \mathcal (\mathcal{U})_j   - \min_{j}  (\mathcal{U})_j \big)^2 \, \mathbb{I},
\end{equation}
with $\mathbb{I} \in \R^{M(M-1) \times M(M-1)}$ the identity matrix. The free parameter $\varsigma_\eta >0$ can be tuned to set the relative noise level. Such a noise model has been used with real-world data,~e.g.,~in~\cite{Darde13a}.

\subsection{Training dataset}\label{sub:traindat}
The training dataset consists of pairs of input vectors and the corresponding output binary values. An input vector of length $M(M-1) = 992$ contains the noisy electrode potentials related to a random choice of the parameter values for the head model presented in Section \ref{ssec:head}. The output is a categorical vector $\{0,1\}$, with label = 1 indicating the presence of a hemorrhage and label = 0 no hemorrhage (there is no stroke or there is an ischemic one). This choice has been motivated by practitioners, since ruling out the presence of a hemorrhage allows them to start treating the patient immediately with blood-thinning medications.

The simulated stroke is represented as a single ball with varying location inside the brain tissue, varying volume and varying conductivity levels drawn from Gaussiann distributions. The chosen parameters $\bar{\sigma}$ and $\varsigma_\sigma$ for all conductivity values are in line with conductivity levels reported in the medical literature \cite{Gabriel96,Lai05,Latikka01,mccann2019variation,Oostendorp00} and are displayed in detail in Table \ref{tab:values}. The radius of the ball defining the stroke is drawn from a uniform distribution $r \sim \mathcal{U}(r_{\rm min}, r_{\rm max})$, with $r_{\rm min} = 0.7$ cm and $r_{\rm max} = 2.3$ cm, which corresponds to volumes ranging from $1.5$ ml to $50$ ml. The inclusion center $(x_c, y_c, z_c )$ is chosen randomly under the condition that the whole inclusion is contained in the brain tissue.

For the azimuthal and polar angles, we have drawn different standard deviations $\varsigma_\theta$ and $\varsigma_\phi$ (cf. Eq~\eqref{eq:Gamma_thph}) uniformly from the set $\{0.01, 0.02, 0.03\}$ (radians) for each forward computation, in order to take into account different levels of electrode movements. These might depend on, e.g., the initial misplacement of the electrode helmet on a patient's head, the differences between the geometry of the patients head and of the helmet, as well as the overall movement of the patient during the examination. For a better understanding on how the selected standard deviation affects the electrode positions, see Figure \ref{fig:electrodes}, where both $\varsigma_\theta$ and $\varsigma_\phi$ are chosen to be $0.03$. In particular, this value is the highest standard deviation one could use in our computational head model before the electrode patches start overlapping, especially on the lower belt, where they are closer.

The relative noise level in \eqref{eq:noisecov} is set to $\varsigma_\eta = 10^{-3}$ (cf.~\cite{Darde13a}, where such a noise level has been used with real-world data). A complete summary of the parameter values used in the random generation of the training data is reported in Table \ref{tab:values}.

The training dataset contains around $40\,000$ samples, approximately split into 50\% conductive inclusions (hemorrhage), 25\% resistive inclusions (ischemia), 25\% no inclusions (healthy). Let us emphasize that all samples, including those associated with no inclusions (strokes) correspond to different head and inner layer shapes, electrode positions and measurement noise realizations. This same statement applies also to all test datasets employed in assessing the performance of the neural networks.

The computations presented were performed with a MATLAB implementation using computer resources within the Aalto University School of Science ``Science-IT" project. Measurements were computed in parallel over $200$ different nodes on Triton \cite{Triton}, the Aalto University high-performance computing cluster, and the overall computation time for generating the training data did not exceed two hours.

\subsection{Test datasets}
While training is performed on the large dataset introduced above, testing the accuracy of the classifiers is realized on independent test sets, based on three different models for the stroke geometries. More precisely, we constructed $14$ test sets that originate from the parameter ranges of Tables~\ref{tab:values} and $4$ different variations, over three families of geometric models for the strokes.

\begin{figure}[]
	\hspace*{\fill}%
	\subcaptionbox{Stroke model (1)}{\includegraphics[width=4.1cm]{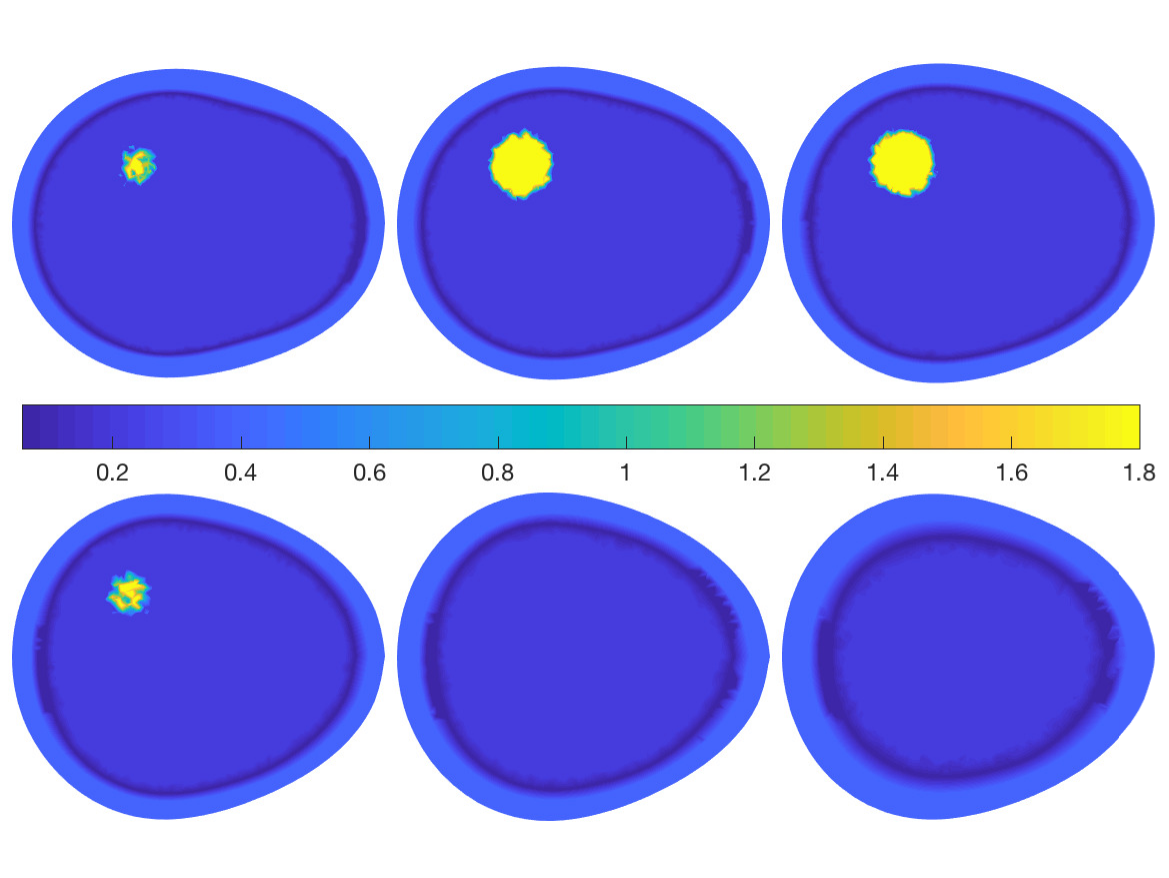}}\hfill%
	\subcaptionbox{Stroke model (2)}{\includegraphics[width=4.1cm]{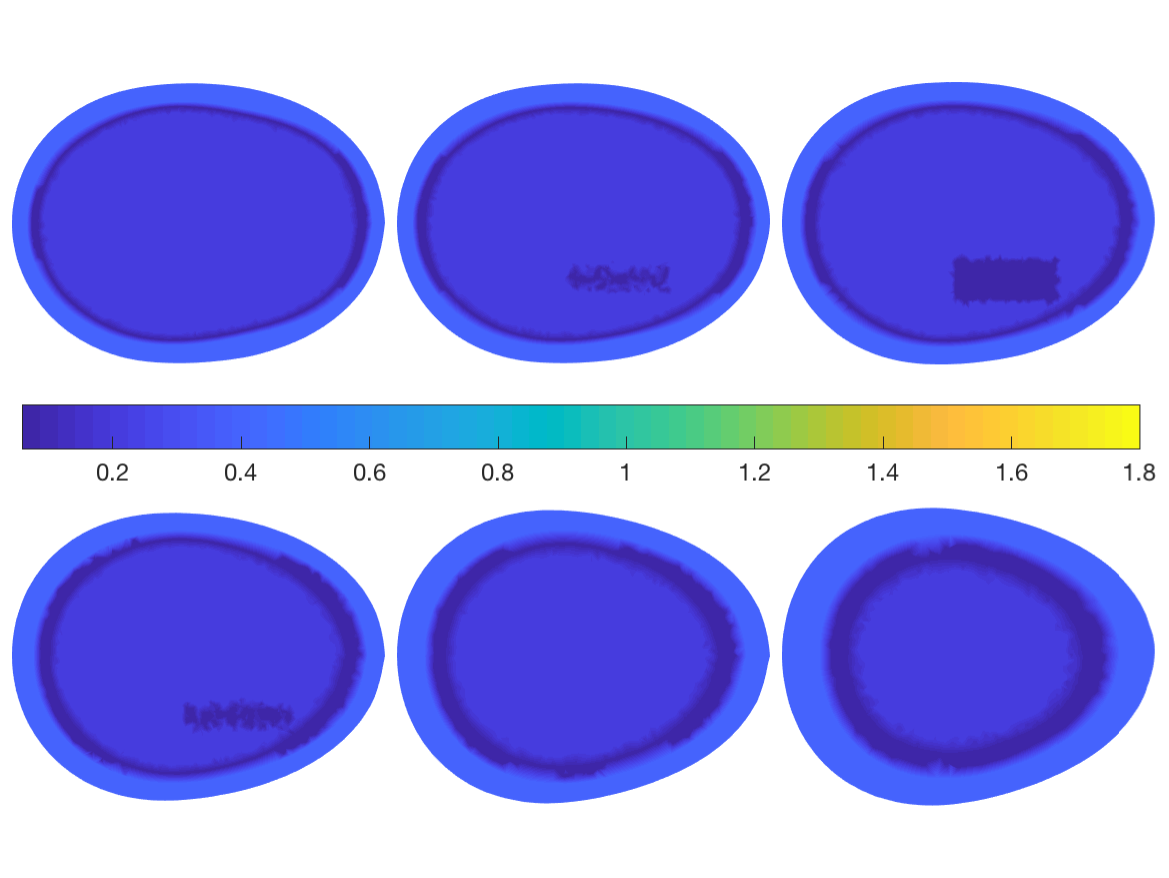}}%
	\hspace*{\fill} \\
	\vspace{0.15cm}
	\hspace*{\fill}%
	\subcaptionbox{Stroke model (3)}{\includegraphics[width=4.1cm]{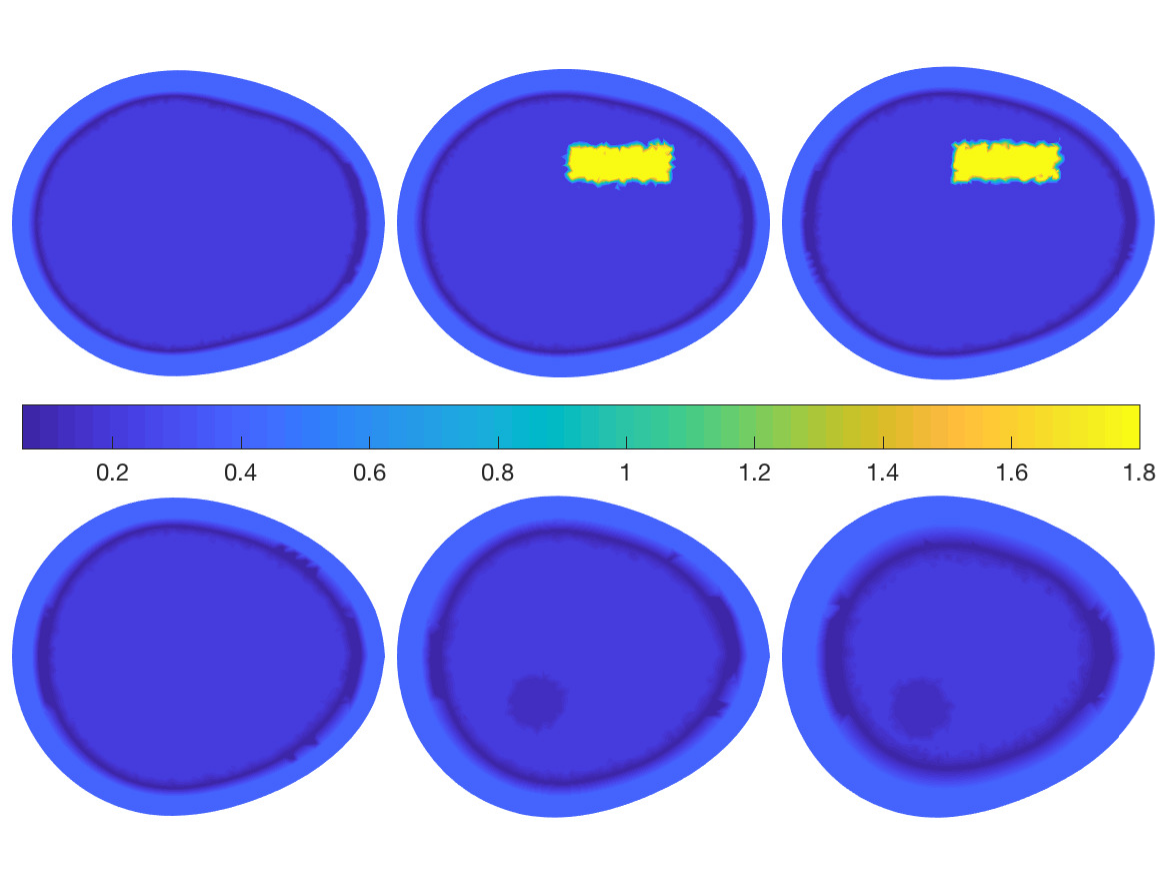}}\hfill%
	\subcaptionbox{Stroke model (3)}{\includegraphics[width=4.1cm]{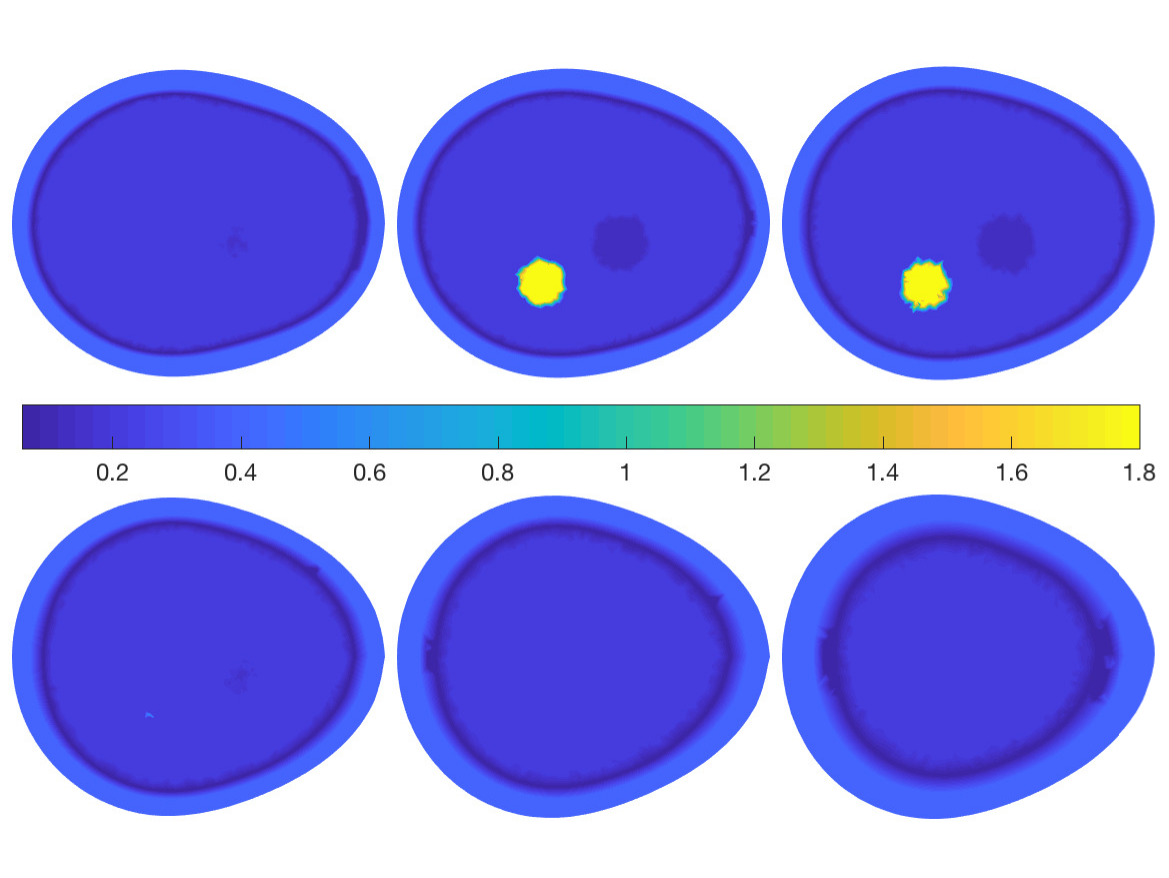}}%
	\hspace*{\fill}%
	\caption{Four examples of six horizontal cross-sections of conductivities used in simulations of data samples in different datasets. The slices are at levels $2$, $3$, $4$, $5$, $6$ and $7$ cm. Top left: a conductivity sample from the training dataset, with one conductive ball (hemorrhagic stroke). Top right: a conductivity sample corresponding to the second geometric stroke model, with one resistive cylinder (ischemic stroke). Bottom row: two examples on the third geometric stroke model, the left panel shows one cylindrical conductive inclusion and one ball-shaped resistive inclusion, while the right panel displays two balls, one conductive and one resistive.}
	\label{fig:sigma}
\end{figure}

The three  models chosen for the test strokes are the following (cf.~Figure \ref{fig:sigma}).
\begin{enumerate}
\item [(1).] A single ball: a test conductivity sample has one ball-shaped inclusion or no inclusion. The corresponding label is equal to $1$ if the inclusion corresponds to a hemorrhagic stroke, $0$ otherwise. We approximately have $50\%$ conductive inclusions, $25\%$ resistive inclusions, $25\%$ no inclusion.
\item [(2).] A single ball or cylinder: a test conductivity sample still exhibits one or no inclusions, but different shapes are considered; a single ball or a cylinder-shaped stroke. The height of the cylinders is uniformly drawn from the distribution $h \sim \mathcal{U}(h_{\rm min}, h_{\rm max})$, with $h_{\rm min} = 1$ cm and $h_{\rm max} = 3$ cm, while the radius is the same as for the corresponding balls, with the values reported in Table \ref{tab:values}. This again corresponds to volumes ranging from $1.5$ ml to $50$ ml. The same labels as in the first case are used. As before, we have $50\%$ conductive inclusions, $25\%$ resistive inclusions and $25\%$ healthy cases; the inclusion shape can be either a ball or a cylinder with a $50\%$ chance.
\item [(3).] Balls or cylinders: we consider cases with one, two or zero inclusions of different shapes (balls and cylinders). The label is equal to $1$, if there is \emph{at least} one hemorrhagic inclusion, and $0$ if the nature of the inclusion(s) is resistive or the brain is healthy. Again, approximately half of the cases have at least one conductive inclusion, $25\%$ of cases have one or two resistive inclusions, $25\%$ have no inclusion. In particular, it is possible that label $1$ corresponds to one hemorrhagic and one resistive inclusion (cf. Figure~\ref{fig:sigma}), which is a case not encountered in the training data.
\end{enumerate}

For each stroke model (1)--(3), we constructed $5$ datasets for which we considered:
\begin{itemize}
\item standard parameters, i.e.,~the ones listed in Table \ref{tab:values},
\item two different ranges for the expected inclusion radius,
\item lower level of electrode misplacement,
\item a higher amount of relative noise added to the data.
\end{itemize}
See Section~\ref{sec:results} for more details. These alterations are to be understood in comparison to the random parameter models listed in Table \ref{tab:values} for the generation of the training data. This results in $14$ different test sets with $5\,000$ samples each. 

\begin{remark}
In this setting, strokes are represented by a well defined ball or cylinder of a constant conductivity value embedded in a homogeneous background, whereas cases of nested inclusions are not considered. In fact, the presence of a penumbra or hypodense tissue would have a huge impact on the performance metrics. A penumbra is a region of normal to high blood volume that surrounds an ischemic stroke as the brain tries to balance the net blood pressure and flow. On the other hand, the hypodense tissue is formed around a hemorrhage when there is shortage of blood corresponding to a lower conductivity. Both situations are critical and detecting them is of vital importance, but we leave their investigation for future studies.
\end{remark}

\section{Results}
\label{sec:results}
We start by first reviewing the computational details about the training and testing process. Next, we present and discuss the classification performance on each dataset.

For each classification learner we evaluate the results in terms of standard \textit{performance metrics} such as sensitivity, specificity and accuracy on the test dataset. Sensitivity, or true positive rate, measures the proportion of actual positives that are correctly identified as such, while specificity, or true negative rate, quantifies the proportion of correctly classified actual negatives. Finally, accuracy is the sum of true positives and true negatives over the total test population, that is, the fraction of correctly classified cases.

\subsection{Training the networks}
The training dataset (see Section~\ref{sub:traindat}) was normalized by subtracting its mean and scaling by its standard deviation before the actual training. The mean and the standard deviation were stored and used to normalize the test datasets in the evaluation phase.

We trained our FCNN and CNN using the Adam optimizer \cite{DBLP:journals/corr/KingmaB14} with batches of size $256$ and learning rate $0.001$. The FCNN was trained for $1500$ epochs on the full training set with no validation. This choice was motivated by preliminary tests with validation that showed that our FCNN was not overfitting the training data. The CNN was trained for $120$ epochs on a $83\% / 17\%$ random training/validation  splitting of the original dataset. In order to show the stability of the classification, the performance metrics for the FCNN and the CNN are reported as the average value computed over $10$ different network trainings.

Training and tests were performed with a Python implementation on a laptop with 8GB RAM and an Intel CPU having clock speed 2.3 GHz. The overall training times were $30$ and $40$ minutes for FCNN and CNN respectively.

\subsection{Results with fully connected neural network}\label{ssec:fnn}
Tables \ref{tab:fnn_D}--\ref{tab:fnn_D2} summarize the results for the trained FCNN and each model for the stroke (1)--(3), respectively. In each table, the highest and lowest accuracy values are highlighted in bold and italics, respectively, while the values concerning the training set in Tables \ref{tab:fnn_D} and \ref{tab:cnn_D} are separated from test cases with a thicker line. Overall, the leading trend is that the network performs considerably well in the case with only one ball-like inclusion, while accuracy somewhat degrades when testing the network in the other two cases, especially for datasets with one or two inclusions of potentially different shapes and type.

The top row of Table \ref{tab:fnn_D} corresponds to testing the network on the training data and thus the resulting accuracy of $0.9673$ is expected to give an upper limit for the performance of the FCNN. However, as shown below, we find that the network performs slightly better on some test datasets. On the top rows of Tables~\ref{tab:fnn_D1} and~\ref{tab:fnn_D2} we present the results for the FCNN on the test sets generated with the same parameters as the training set, but with different geometric models for the stroke. In every table, the results on the other rows are obtained by testing the trained network on four different test sets, where each time some parameters are altered from the values in Table \ref{tab:values}. 

The second rows of all three tables correspond to testing our classifier with strokes of larger volumes, where the lower bound for the ball radius is increased to $r_{\rm min} = 1$ cm, leading to the inclusion radius being drawn from $r \sim \mathcal{U}(1, 2.3)$ cm (4.20 ml -- 50 ml for ball-shaped inclusions). In the case of a cylindrical inclusion, its radius is modified accordingly, while its height remains in the range $[1,3]$ cm, as in Table \ref{tab:values} (9.40 ml -- 50 ml for cylinder-shaped inclusions). This choice of parameters leads to a better performance of the network, due to the better average visibility of the test inclusions to EIT measurements. Accuracy ranges from $91\%$ in the case (3) with one or more inclusions of different shapes up to $97\%$ accuracy for the test (1) with only one ball inclusion, which is actually higher than the accuracy on the training set.

Conversely, as shown on the third rows of each table, when only smaller inclusions are considered, the overall performance metrics decrease, with the lowest accuracy of $87\%$ reached for the last geometric stroke model (3). Again, the height of the possible cylinder-shaped inclusions remains as in Table \ref{tab:values}. It should also be noted that in Table \ref{tab:fnn_D}, as well as in all our other tests, test datasets with larger and smaller inclusions are not simply subsets of the training set, but they have instead been simulated with random choices for the remaining parameters (cf.~Table \ref{tab:values}), providing hence new and unseen sets of data for the neural network.

In the experiments documented on the fourth rows, we considered a lower degree of electrode misplacement, that is, a lower level of inaccuracy in the electrode position when a helmet of electrodes is placed on a patient's head. This translates to selecting a lower standard deviation in formula \eqref{eq:Gamma_thph} when simulating the test data, with $\varsigma_\theta, \varsigma_\phi $ being drawn from $ \{0.005, 0.01, 0.02\}$ radians. Results in this case are notably better than those obtained with the baseline dataset for each of the three random stroke models (1)--(3), with the accuracy ranging from $91\%$ to almost $97\%$. This confirms that electrode positioning does indeed significantly affect the classification accuracy. 

Finally, to test the network in a more realistic setup, we increased the relative noise level in the measurements in \eqref{eq:noisecov} from $\varsigma_\eta = 10^{-3}$ to $\varsigma_\eta = 10^{-2}$ in the test data. The resulting performance indicators are listed on the last rows of Tables \ref{tab:fnn_D}--\ref{tab:fnn_D2}. Despite the (inevitable) decrease in performance, the accuracy was still over $90\%$ for all three models for test stroke generation (1)--(3). 

%DATASET D
\begin{table}[]
\caption{Means and standard deviations for the performance metrics obtained after training our FCNN $10$ times on the same $40\,000$-sample training dataset and testing it on different test sets of $5\,000$ samples with a single ball-shaped inclusion. Note that the first row corresponds to testing the FCNN on the training dataset.}
\begin{tabular}{lllllll}
\hline
\multicolumn{7}{c}{FCNN: Datasets with a single ball-shaped stroke}\\ \hline
\multirow{2}{*}{Datasets}  & \multicolumn{2}{c}{Accuracy} & \multicolumn{2}{c}{Sensitivity} & \multicolumn{2}{c}{Specificity} \\ \cline{2-7} 
& Mean & \multicolumn{1}{c}{\begin{tabular}[c]{@{}c@{}} Std\\ Dev\end{tabular}} & Mean & \multicolumn{1}{c}{\begin{tabular}[c]{@{}c@{}}Std\\ Dev\end{tabular}} & Mean & \multicolumn{1}{c}{\begin{tabular}[c]{@{}c@{}}Std\\ Dev\end{tabular}} \\ \hline
Training dataset   &   0.9673   &  0.0034 &   0.9753   &  0.0116   &   0.9600  &  0.0122  \\ \noalign{\hrule height 1.3pt}
Larger radius $r\in [1, 2.3] $ cm &   \textbf{0.9740}   & 0.0040   &   0.9660   &  0.0107 &   0.9825   &  0.0054  \\ %\hline
Smaller radius  $r\in [0.7, 1.5] $ cm &  \textit{0.9190}  & 0.0081  &   0.9500   &  0.0181    &   0.8933   &   0.0243 \\ %\hline
Less el. uncert. $\varsigma_\theta, \varsigma_\phi \in \{0.5, 1, 2\}\cdot 10^{-2}$rad &   0.9684   &  0.0039  &   0.9769   &  0.0121  &    0.9609  &   0.0107  \\ %\hline
Increased relative noise $\varsigma_\eta = 10^{-2}$    &  0.9527   &   0.0033   &   0.9617  &   0.0129   &   0.9445   &   0.0124 \\ \hline
\end{tabular}
\label{tab:fnn_D}
\end{table}

%DATASET D1
\begin{table}[]
\caption{Performance metrics after testing the FCNN on test sets with a single ball-shaped or cylinder-shaped inclusion. Compared to the values in Table \ref{tab:fnn_D}, accuracy slightly drops but the trend remains the same, with tests on dataset corresponding to inclusions with larger radii being the most accurate and the ones with smaller radii displaying only a $89\%$ accuracy.}
\begin{tabular}{lllllll}
\hline
\multicolumn{7}{c}{FCNN: Datasets with a single ball or cylinder-shaped stroke}\\ \hline
\multirow{2}{*}{Datasets}  & \multicolumn{2}{c}{Accuracy} & \multicolumn{2}{c}{Sensitivity} & \multicolumn{2}{c}{Specificity} \\ \cline{2-7} 
& Mean & \multicolumn{1}{c}{\begin{tabular}[c]{@{}c@{}}Std\\ Dev\end{tabular}} & Mean & \multicolumn{1}{c}{\begin{tabular}[c]{@{}c@{}}Std\\ Dev\end{tabular}} & Mean & \multicolumn{1}{c}{\begin{tabular}[c]{@{}c@{}}Std\\ Dev\end{tabular}} \\ \hline
Standard parameters   &   0.9347   & 0.0074  &   0.9460  &  0.0093   &   0.9243   &  0.0148 \\ %\hline
Larger radius $r\in [1, 2.3] $ cm &   \textbf{0.9523}   &  0.0062  &  0.9439    & 0.0101  &   0.9609   &  0.0116  \\ %\hline
Smaller radius  $r\in [0.7, 1.5] $ cm &   \textit{0.8988}   &  0.0128   &  0.9398    &   0.0110   &    0.8644  &  0.0240  \\ %\hline
Less el. uncert. $\varsigma_\theta, \varsigma_\phi \in \{0.5, 1, 2\}\cdot 10^{-2}$rad &   0.9412   &  0.0066  &   0.9508   &  0.0097  &   0.9323   &  0.0142   \\ %\hline
Increased relative noise $\varsigma_\eta = 10^{-2}$   &   0.9335   &   0.0081   &  0.9454    &   0.0100   &   0.9231   &   0.0155 \\ \hline
\end{tabular}
\label{tab:fnn_D1}
\end{table}

%DATASET D2
\begin{table}[]
\caption{Performance metrics after testing the FCNN on test sets with ball-shaped or cylinder-shaped inclusions. Levels of accuracy are lower compared to Tables \ref{tab:fnn_D} and \ref{tab:fnn_D1} but, once again, the overall trend remains unchanged.}
\begin{tabular}{lllllll}
\hline
\multicolumn{7}{c}{FCNN: Datasets with balls or cylinders-shaped strokes}\\ \hline
\multirow{2}{*}{Datasets}  & \multicolumn{2}{c}{Accuracy} & \multicolumn{2}{c}{Sensitivity} & \multicolumn{2}{c}{Specificity} \\ \cline{2-7} 
& Mean & \multicolumn{1}{c}{\begin{tabular}[c]{@{}c@{}}Std\\ Dev\end{tabular}} & Mean & \multicolumn{1}{c}{\begin{tabular}[c]{@{}c@{}}Std\\ Dev\end{tabular}} & Mean & \multicolumn{1}{c}{\begin{tabular}[c]{@{}c@{}}Std\\ Dev\end{tabular}} \\ \hline
Standard parameters   &  0.9038    & 0.0088  & 0.8954    &   0.0210  &  0.9142    &  0.0182 \\ %\hline
Larger radius $r\in [1, 2.3] $ cm  &   0.9135   &  0.0093  &  0.8945    &  0.0204 &   0.9351   &  0.0151  \\ %\hline
Smaller radius  $r\in [0.7, 1.5] $ cm &  \textit{0.8793}    &   0.0108  &   0.8867   &   0.0237   &   0.8740   &  0.0205  \\ %\hline
Less el. uncert. $\varsigma_\theta, \varsigma_\phi \in \{0.5, 1, 2\}\cdot 10^{-2}$rad &  \textbf{0.9166}   & 0.0102   &   0.9055   &  0.0235  &   0.9297  &   0.0145  \\ %\hline
Increased relative noise $\varsigma_\eta = 10^{-2}$   &   0.9020   &  0.0094    &  0.8933  &   0.0200   & 0.9127    &   0.0171 \\ \hline
\end{tabular}
\label{tab:fnn_D2}
\end{table}

\subsection{Results with convolutional neural network}\label{ssec:cnn}
Following the same workflow as in Section \ref{ssec:fnn}, Tables \ref{tab:cnn_D}--\ref{tab:cnn_D2} present the analogous results for the CNN trained on the training dataset and tested in the three cases related to the different geometric models for the inclusions. It can clearly be seen that the overall accuracy is significantly lower than in the case of FCNN. This is arguably due to the fact that our CNNs tend to overfit the data, resulting in a high accuracy for the training set (first row of Table \ref{tab:cnn_D}), which significantly degrades in the other test cases. With the exception of the case of the smaller strokes, the overall accuracy remains steadily over $80\%$ and roughly follows the same trends as for the FCNN. Tests in the case of larger inclusion volumes have by far the best accuracy level, ranging from $84\%$ to almost $89\%$, while the least accurate classification is indeed obtained for the strokes with smaller volumes. The mean accuracy with a lower level of electrode uncertainty and increased relative noise is between $81\%$ and $86\%$, depending on the model for generating the test strokes.

%DATASET D
\begin{table}[H]
\caption{Means and standard deviations for the performance metrics obtained after training our CNN $10$ times on the same $40\,000$-sample training dataset and testing it on different test sets of $5\,000$ samples with a single ball-shaped inclusion. Note that the first row corresponds to testing the CNN on the training dataset.}
\begin{tabular}{lllllll}
\hline
\multicolumn{7}{c}{CNN: Datasets with a single ball-shaped stroke}\\ \hline
\multirow{2}{*}{Datasets}  & \multicolumn{2}{c}{Accuracy} & \multicolumn{2}{c}{Sensitivity} & \multicolumn{2}{c}{Specificity}  \\ \cline{2-7} 
& Mean & \multicolumn{1}{c}{\begin{tabular}[c]{@{}c@{}}Std\\ Dev\end{tabular}} & Mean & \multicolumn{1}{c}{\begin{tabular}[c]{@{}c@{}}Std\\ Dev\end{tabular}} & Mean & \multicolumn{1}{c}{\begin{tabular}[c]{@{}c@{}}Std\\ Dev\end{tabular}} \\ \hline
Training dataset  &   \textbf{0.9593}   & 0.0148  &   0.9752   &  0.0215   &   0.9479  &  0.0383 \\ \noalign{\hrule height 1.3pt}
Larger radius $r\in [1, 2.3] $ cm &  0.8892   &  0.0096  &   0.8747   & 0.0399  &  0.9110   &  0.0278  \\ %\hline
Smaller radius  $r\in [0.7, 1.5] $ cm  &  \textit{0.7589}    &  0.0237   &    0.8435  &  0.0401    &   0.7152   &  0.0452  \\ %\hline
Less el. uncert. $\varsigma_\theta, \varsigma_\phi \in \{0.5, 1, 2\}\cdot 10^{-2}$rad &   0.8611  &  0.0040  &   0.8823 &  0.0398  & 0.8489 &   0.0329  \\ %\hline
Increased relative noise $\varsigma_\eta = 10^{-2}$  &  0.8448  &  0.0041    &  0.8664  &  0.0403    &  0.8324  &    0.0338 \\ \hline
\end{tabular}
\label{tab:cnn_D}
\end{table}

%DATASET D1
\begin{table}[H]
\caption{Performance metrics after testing the CNN on test sets with a single ball-shaped or cylinder-shaped inclusion. Compared to the values in Table \ref{tab:cnn_D}, accuracy remains somewhat stable, with tests on the dataset corresponding to inclusions with larger radii being the most accurate and the ones with smaller radii displaying only a $76\%$ accuracy.}
\begin{tabular}{lllllll}
\hline
\multicolumn{7}{c}{CNN: Datasets with a single ball or cylinder-shaped stroke}\\ \hline
\multirow{2}{*}{Datasets}  & \multicolumn{2}{c}{Accuracy} & \multicolumn{2}{c}{Sensitivity} & \multicolumn{2}{c}{Specificity} \\ \cline{2-7} 
& Mean & \multicolumn{1}{c}{\begin{tabular}[c]{@{}c@{}}Std\\ Dev\end{tabular}} & Mean & \multicolumn{1}{c}{\begin{tabular}[c]{@{}c@{}}Std\\ Dev\end{tabular}} & Mean & \multicolumn{1}{c}{\begin{tabular}[c]{@{}c@{}}Std\\ Dev\end{tabular}} \\ \hline
Standard parameters   &  0.8470    &  0.0057 &   0.8677    &   0.0299  &   0.8328  &  0.0294 \\ %\hline
Larger radius $r\in [1, 2.3] $ cm &    \textbf{0.8756}  & 0.0045   &  0.8660    & 0.0322  & 0.8902   &  0.0300  \\ %\hline
Smaller radius  $r\in [0.7, 1.5] $ cm &  \textit{0.7608}    &  0.0234   &  0.8465    &  0.0275    &   0.7120   &   0.0405 \\ %\hline
Less el. uncert. $\varsigma_\theta, \varsigma_\phi \in \{0.5, 1, 2\}\cdot 10^{-2}$rad &  0.8508  &  0.0049  & 0.8776  & 0.0280   & 0.8312  &  0.0282   \\ %\hline
Increased relative noise $\varsigma_\eta = 10^{-2}$  &  0.8477   &  0.0034    &  0.8646  &   0.0322   &  0.8378  & 0.0289    \\ \hline
\end{tabular}
\label{tab:cnn_D1}
\end{table}

%DATASET D2
\begin{table}[H]
\caption{Performance metrics after testing the CNN on test sets with ball-shaped or cylinder-shaped inclusions. Levels of accuracy are lower compared to Tables \ref{tab:cnn_D} and \ref{tab:cnn_D1} but, with the exception of the case of smaller strokes, the overall accuracy remains steadily over $80\%$.}
\begin{tabular}{lllllll}
\hline
\multicolumn{7}{c}{CNN: Datasets with balls or cylinders-shaped strokes}\\ \hline
\multirow{2}{*}{Datasets}  & \multicolumn{2}{c}{Accuracy} & \multicolumn{2}{c}{Sensitivity} & \multicolumn{2}{c}{Specificity} \\ \cline{2-7} 
& Mean & \multicolumn{1}{c}{\begin{tabular}[c]{@{}c@{}}Std\\ Dev\end{tabular}} & Mean & \multicolumn{1}{c}{\begin{tabular}[c]{@{}c@{}}Std\\ Dev\end{tabular}} & Mean & \multicolumn{1}{c}{\begin{tabular}[c]{@{}c@{}}Std\\ Dev\end{tabular}} \\ \hline
Standard parameters   &   0.8181    & 0.0247  &   0.7826    &  0.0504   &  0.8750  &  0.0203 \\ %\hline
Larger radius $r\in [1, 2.3] $ cm &   \textbf{0.8446}   & 0.0304   &   0.7916  &  0.0514 &   0.9292   &  0.0163  \\ %\hline
Smaller radius  $r\in [0.7, 1.5] $ cm  &  \textit{0.7510}    &  0.0103   &  0.7548    &   0.0486   &  0.7590    &   0.0310 \\ %\hline
Less el. uncert. $\varsigma_\theta, \varsigma_\phi \in \{0.5, 1, 2\}\cdot 10^{-2}$rad &  0.8362  & 0.0238   & 0.7991  &  0.0483  &  0.8918  &  0.0178   \\ %\hline
Increased relative noise $\varsigma_\eta = 10^{-2}$  &   0.8144  &   0.0261   &  0.7840  &   0.0523   &  0.8646  & 0.0192    \\ \hline
\end{tabular}
\label{tab:cnn_D2}
\end{table}

\section{Conclusions}
\label{sec:conclusion}
This work applies neural networks to the detection of brain hemorrhages from simulated absolute EIT data on a $3$D head model. We developed large datasets based on realistic anatomies that we used to train and test a fully connected neural network and a convolutional neural network. Our classification tests show encouraging results for further development of these techniques, with fully connected neural networks achieving an average accuracy higher than $90\%$ on most test datasets. Since the datasets included a large proportion of healthy patents, this work could motivate the development of a new, cheap, non-invasive screening test for brain hemorrhages.

The results demonstrate that the classification performance is affected by many factors, including the size of the stroke, the mismodeling of the electrode position and the noise in the data. Among these, the size of the stroke is the one most significantly altering the performance, which indicates that the method is probably unable to detect very small hemorrhages. Note that in the dataset generation we considered strokes with volumes of as small as 1.5 ml. Further, the FCNN, trained on single ball-shaped strokes, was able to generalize to data generated from multiple strokes of different shapes with little loss in the accuracy. We also want to stress that in every dataset the position of the stroke was chosen uniformly at random within the brain tissue. This means that on the one hand we did not restrict or assume to know in which brain hemisphere the stroke was located. On the other hand we considered strokes taking place potentially very close to the skull layer or very deep inside the brain tissue, a situation which is very challenging to analyze from EIT data.

The simulated datasets have several limitations nonetheless. The three-layer head model is a simplified version of the true head anatomy. In particular, it does not take into account the cerebrospinal fluid and other subtle anatomical tissues. Besides, the considered conductivity distributions have an average skull-to-brain ratio of $3:10$, which has been shown to be a difficult quantity to estimate \cite{Lai05,Oostendorp00} and appears to affect the performance metrics. A brain stroke is a complex phenomenon which evolves with time, and it is so far unclear how to precisely model its conductivity values, both in the hemorrhagic and in the ischemic case.

This work leaves many open directions for future studies. A natural one is to consider more realistic datasets. To our knowledge, the only publicly available dataset with real EIT data on human patients has been released by University College London \cite{goren2018multi}. The dataset comprises data from only 18 patients, which is far from the size of our training set (40\,000 samples), but could be used as a test set after training on synthetic data. Datasets from phantom data would also be an important intermediate step to validate the generation of synthetic measurements. On a related note, cases of nested inclusions should also be taken into account.

Another important direction is to implement more refined machine learning algorithms. We found that a simple 2-layer FCNN ourperforms a more sophisticated CNN, though we did not investigate in depth the problem of designing an optimal network architecture for our EIT data. This certainly leaves room for improvements in the performance.

Finally, bedside real-time monitoring of stroke patients could become another valuable application of the presented techniques. The aim would be to predict whether an hemorrhage is increasing in volume over time or it is stable. This could be studied with a similar approach based on machine learning algorithms trained on datasets with several measurements taken at different times.

\section*{Acknowledgments}
This work has been partially carried out at the Machine Learning Genoa (MaLGa) center, Universit\`a di Genova (IT). This material is based upon work supported by the Air Force Office of Scientific Research under award number FA8655-20-1-7027. MS is member of the ``Gruppo Nazionale per l'Analisi Matematica, la Probabilit\`a e le loro Applicazioni'' (GNAMPA), of the ``Istituto Nazionale per l'Alta Matematica'' (INdAM). This work was also supported by the Academy of Finland (decision 312124) and the Finnish Cultural Foundation (grant number 00200214). The authors thank Nuutti Hyv\"onen (Aalto University) for useful discussions, Antti Hannukainen (Aalto University) for help with his FE solver and all the members of the project ``Stroke classification and monitoring using Electrical Impedance Tomography'', funded by Jane and Aatos Erkko Foundation, that motivated this work.

\section*{Conflict of interest}
All authors declare no conflicts of interest in this paper.

%\section*{Supplementary (if necessary)}

\end{document}